\documentclass[reqno,10pt,a4paper]{amsart}
\usepackage{amsmath}
\usepackage{amssymb}
\usepackage{latexsym}
\usepackage{epsf}

\usepackage[all]{xy}
\UseComputerModernTips
\CompileMatrices
\newcommand{\xyinto}{\ar@{ >->}}
\newcommand{\xyonto}{\ar@{>>}}
\newcommand{\xyequals}{\ar@{=}}

\DeclareMathOperator{\Hom}{Hom}
\DeclareMathOperator{\Ext}{Ext}
\DeclareMathOperator{\End}{End}
\renewcommand{\ge}{\geqslant}
\renewcommand{\le}{\leqslant}

\newcommand{\N}{\mathbb{N}}

\DeclareMathOperator{\glob}{glob}
\newcommand{\frob}{^{\mathrm F}}
\newcommand{\rarr}{\rightarrow}

\DeclareMathOperator{\inj}{inj}

\newcommand{\B}{{\mathcal B}}
\DeclareMathOperator{\good}{{\mathcal F}(\nabla)}
\DeclareMathOperator{\doog}{{\mathcal F}(\Delta)}

\DeclareMathOperator{\gfd}{gfd}
\DeclareMathOperator{\wfd}{wfd}
\newcommand{\impl}{\Rightarrow}
\newcommand{\lrimpl}{\Leftrightarrow}
\DeclareMathOperator{\hd}{hd}
\DeclareMathOperator{\soc}{soc}
\newcommand{\Mod}{\mathrm{mod}}
\DeclareMathOperator{\hw}{hw}
\DeclareMathOperator{\Ind}{Ind}
\DeclareMathOperator{\partn}{\Lambda^+}
\newcommand{\GL}{\mathrm{GL}}
\newcommand{\SL}{\mathrm{SL}}
\newcommand{\St}{\mathrm{St}}

\newcommand{\Zone}{\hat{Z}^\prime_1}
\newcommand{\Lone}{\hat{L}_1}
\newcommand{\uone}{\mathrm{u}_1}

\begin{document}
\theoremstyle{plain}
\newtheorem{thm}{Theorem}[section]
\newtheorem{propn}[thm]{Proposition}
\newtheorem{cor}[thm]{Corollary}
\newtheorem{clm}[thm]{Claim}
\newtheorem{lem}[thm]{Lemma}
\newtheorem{conj}[thm]{Conjecture}
\theoremstyle{definition}
\newtheorem{defn}[thm]{Definition}
\newtheorem{rem}[thm]{Remark}

\title{The Global Dimension of Schur Algebras for $\GL_2$ and $\GL_3$}
\author{Alison E. Parker}
\address{School of Mathematical Sciences, Queen Mary 
College, Mile End Road, London E1 4NS, UK}
\email{A.Parker@qmw.ac.uk}

\begin{abstract}{We first define the notion of good filtration
dimension and Weyl filtration dimension in a quasi-hereditary
algebra.
We calculate these dimensions explicitly for all irreducible modules
in $\SL_2$ and $\SL_3$.
We use these to show that the global dimension
of a Schur algebra for $\GL_2$  and $\GL_3$ is twice 
the good filtration dimension. 
To do this for $\SL_3$, we give an explicit filtration of
the modules $\nabla(\lambda)$ by modules of the form $\nabla(\mu)\frob
\otimes L(\nu)$ where $\mu$ is a dominant weight and $\nu$ is
$p$-restricted. }
\end{abstract}


\maketitle

\section*{Introduction and Background}

The global dimension of a $q$-Schur algebra $S_q(n,r)$ has been determined
when $r\le n$. 
This was calculated by Totaro~\cite{totaro} (for the
classical case) and Donkin~\cite{donkbk} (for the quantum case).
In general the global dimension of $S_q(n,r)$ is not known, although
we do have upper bounds for this value, (see~\cite{totaro} for more
details). 
In this paper we will calculate explicitly the global dimension of
$S(n,r)$ for $n=2$ and $n=3$. We also find the global dimension of
$S_q(2,r)$.

We first briefly review some of the notation and definitions that we
will use in this paper. 
The reader is referred to~\cite{donkbk} and~\cite{jantz} for further information. 
We let $G=\SL_n(k)$ where $k$ is an algebraically closed field of
characteristic $p$ and $\mathrm{F}:G \rarr G$ the corresponding
Frobenius morphism. We let $G_1$ be the first Frobenius
kernel.
Let $T$ be a maximal torus of $G$, $W$ the corresponding Weyl group
and $B \supseteq T$ a Borel subgroup. Let $X= X(T)$ be the weight
lattice and let $X^+$ be the set of dominant weights. 
Let $X_1$ be the set of $p$-restricted dominant weights
and $A_0$ the set of weights in the interior of the fundamental
alcove.

For $\lambda \in X^+$, let $k_\lambda$ be the one-dimensional module
for $B$ which has weight $\lambda$. We define $\nabla(\lambda)=
\Ind_B^G(k_\lambda)$.
This module has character given by Weyl's character formula and has
simple socle $L(\lambda)$, the irreducible $G$-module of highest weight
$\lambda$.

We define $S(n,r)$ to be the Schur algebra corresponding to $\GL_n$
and $\partn(n,r)$ to be the set of partitions of $r$ with
at most $n$ parts.  
Modules for $S(n,r)$ are naturally polynomial modules for $\GL_n$ which are
homogeneous of degree $r$. An irreducible module for $S(n,r)$ 
is also an irreducible module for $\GL_n$ and correspond to elements
of $\partn(n,r)$. We write $L(a_1,\cdots,a_n)$ for the irreducible
module of highest weight $(a_1,\cdots,a_n)\in \partn(n,r)$.
More details can be found in~\cite[chapter 2]{green}. 
We have a natural correspondence between weights for $S(n,r)$-modules
and $\SL_n$-modules given by 
$$ \lambda= (\lambda_1,\lambda_2, \ldots, \lambda_n) \mapsto
(\lambda_1-\lambda_2, \lambda_2-\lambda_3, \ldots,
\lambda_{n-1}-\lambda_n).$$
The usual partial ordering on $\SL_n$-weights is equivalent by the above
correspondence to the dominance ordering on partitions.

The category of rational $G$-modules has enough injectives and so we may
define $\Ext^*(-,-)$ as usual by using injective resolutions.
We have a canonical isomorphism for all
$S(n,r)$-modules~\cite[2.2d]{donk1}
$$ \Ext_{S(n,r)}^i(V,W) \cong\Ext^i_{\GL_n}(V,W) \cong \Ext^i_{\SL_n}(V,W). $$
Thus we may do all $\Ext$ calculations in $\Mod(G)$.


\section{Quasi-hereditary Algebras}
We start with the definition of a highest weight category 
given by Cline, Parshall and Scott~\cite{CPS}.
\begin{defn}
Let $k$ be a field, $S$ a finite dimensional algebra over $k$,
$\Lambda$ an indexing set for the isomorphism classes of simple
$S$-modules
with a correspondence $\lambda \leftrightarrow L(\lambda)$,
and $\le$ a partial order on $\Lambda$. We say $(S, \Lambda)$ is a
\emph{highest weight category} if and only if, for all $\lambda \in
\Lambda$ there is a left $S$-module $\nabla(\lambda)$, called the
\emph{costandard module}, such that:
\begin{enumerate}
\item[(i)]there exists an injection $\phi_\lambda: L(\lambda)
\rarr\nabla(\lambda)$, and the composition factors, $L(\mu)$ of the
cokernel satisfy $\mu < \lambda$ 
\item[(ii)]the indecomposable injective hull, $I(\lambda)$, of
$L(\lambda)$ contains $\nabla(\lambda)$ via the injection $\psi_\lambda:
\nabla(\lambda) \rarr I(\lambda)$ and the cokernel of $\psi_\lambda$ is
filtered by modules $\nabla(\mu)$ with $\mu > \lambda$.
\end{enumerate}
\end{defn}
Dually we have $\Delta(\lambda)$ as standard modules by replacing
injective by projective and cokernel by kernel in the above definition.
Thus the standard modules  $\Delta(\lambda)$ have simple head
$L(\lambda)$.
(We will assume that $S$ is Schurian in the sense that $\End_S(L)=k$
for all simple modules $L$.)

A highest weight category is equivalent to the module category
for a quasi-hereditary
algebra although we will not
show this here. See~\cite{klukoe} or~\cite[Appendix]{donkbk} for
details.

We say $X \in \Mod(S)$ has a \emph{good filtration} if it
has a filtration
$$0=X_0 \subset X_1 \subset \cdots \subset X_i=X$$ 
with quotients $X_j/X_{j-1}$ isomorphic to 
$\nabla(\mu_j)$ for some $\mu_j \in \Lambda$. We denote the class of
$S$-modules with good filtration $\good$, and dually the class of
modules filtered by $\Delta(\mu)$'s as $\doog$. We will say that $X
\in \doog$ has a \emph{Weyl filtration}. The multiplicity of
$\nabla(\mu)$ in a filtration of $X \in \good$ is independent of the
filtration chosen. This multiplicity is 
denoted by $\bigl(X:\nabla(\mu)\bigr)$ and  the composition
multiplicity of $L(\mu)$ in $ X\in \Mod(S)$ is denoted by $\big[X:L(\mu)\big]$.

We state some of the results that we will need
from~\cite[Appendix A]{donkbk}.

\begin{propn}
\begin{enumerate}
\item[(i)]\
Let $X \in \Mod(S)$ and $\lambda \in \Lambda$. If $\Ext^1\bigl(X,
\nabla(\lambda)\bigr)\ne 0$ then $X$ has a composition factor $L(\mu)$ with
$\mu > \lambda$.
\item[(ii)]\label{propn:good}\
For $X \in \doog$, $Y\in \good$ and $i>0$, we have 
$\Ext^i(X,Y)=0$.
\item[(iii)]\
Suppose $\Ext^1\bigl(\Delta(\mu),M\bigr)=0$ for all $\nu \in \Lambda$
then $M\in \good$.
\item[(iv)]\
Let $X \in \good$ (resp. $X \in \doog $) and $Y$ a direct summand of
$X$ then $Y\in \good$ (resp. $Y \in \doog$).
\end{enumerate}
\end{propn}
\begin{proof}
See
~\cite[A2.2]{donkbk}.
\end{proof}


%

\section{Good filtration, Weyl and global dimensions}

In this section $S$ is a quasi-hereditary algebra with poset
$(\Lambda,\le)$.


Suppose $X \in \Mod (S)$. We can resolve $X$ by modules $M_i\in \good$ as
follows
$$0\rarr X \rarr M_0 \rarr M_1 \rarr \dots \rarr M_d \rarr 0.$$
We call such a resolution a \emph{good resolution} for $X$.
Good resolutions exist for all $S$-modules as an injective
resolution is also a good resolution.

\begin{defn}
Let $X \in \Mod (S)$. We say $X$ has \emph{good filtration dimension}
$d$, denoted  $\gfd(X)=d$,
if the following two equivalent conditions hold:
\begin{enumerate}
\item[(i)]
$0\rarr X \rarr M_0 \rarr M_1 \rarr \dots \rarr M_d \rarr 0$
is a resolution for $X$ with $M_i \in \good$, of shortest
possible length.
\item[(ii)]
$\Ext^i(\Delta(\lambda),X)=0$ for all $i>d$ and all $\lambda \in
\Lambda$,
but there exists $\lambda \in \Lambda$ such that
$\Ext^d(\Delta(\lambda),X)\ne 0$.
\end{enumerate}
\end{defn}

\begin{proof} 
See~\cite[proposition 3.4]{fripar}.
\end{proof}

Similarly we have the dual notion of the 
\emph{Weyl filtration dimension} of $M$ which we will denote 
$\wfd(M)$. 

\begin{lem}\label{lem:wandg}
Given $S$-modules $M$ and $N$, we have
$$\Ext^i(N,M)=0\mbox{ for }i > \wfd(N)+\gfd(M).$$
\end{lem}
\begin{proof} We proceed by induction on $\wfd(N)+\gfd(M)$.
If either of $\wfd(N)=0$ or $\gfd(M)=0$
then we are done by the definition of good filtration dimension and
Weyl filtration dimension, so we may assume that both $\wfd(N)$ and
$\gfd(M)$ are non-zero.
Now we can embed $M$ in $A\in\good$ with
quotient $B$.  
The exact sequence gives us:
$$\cdots \rarr \Ext^{i-1}(N,B)\rarr \Ext^i(N,M) \rarr \Ext^i(N,A)
\rarr \cdots .$$
Now $\Ext^i(N,A) = 0 $ if $i> \wfd(N).$ Also $\Ext^{i-1}(N,B)=0$ as
$i-1>\wfd(N)+\gfd(B)$ by induction. 
($B$ has strictly smaller good filtration dimension
than $M$). 
Hence $\Ext^i(N,M)=0$. 
\end{proof}

\begin{defn}
Let $g=\sup\{ \gfd(X) \mid X \in \Mod(S)\}$. 
We say $S$ has good filtration
dimension $g$ and denote this by $\gfd(S)=g$.
Let $w=\sup\{ \wfd(X) \mid X \in \Mod(S)\}$. 
We say $S$ has Weyl filtration
dimension $w$ and denote this by $\wfd(S)=w$.
\end{defn}

\begin{rem} In general $\gfd(S)$ is not the good filtration dimension
of $S$ when considered as its own left (or right) module. Similar
remarks apply to $\wfd(S)$. We will only
use $\gfd(S)$ and $\wfd(S)$ in the sense that they are defined above.
\end{rem}

For a finite dimensional $k$-algebra $S$, 
the \emph{injective dimension} of an $S$-module $M$, 
is the length of a shortest possible injective
resolution and is denoted by $\inj(M)$.
Equivalently we have
$$\inj(M)=\sup\{ d \mid \Ext^d(N,M)\not\cong 0 \mbox{ for } N \in \Mod(S)\}.$$
 The \emph{global dimension} of $S$ is the
supremum of all the injective dimensions for $S$-modules, and is
denoted by $\glob(S)$. This is equivalent to
$$\glob(S)=\sup\{ d \mid \Ext^d(N,M) \not\cong 0 \mbox{ for } N,M \in \Mod(S)\}.$$

\begin{cor}\label{cor:glob}
The global dimension of
$S$ has an upper bound of $\wfd(S)+\gfd(S)$.
\end{cor}

Now suppose $S$ is a quasi-hereditary algebra with contravariant
duality preserving simples.
That is there exists an involutory, contravariant functor 
$^\circ: \Mod (S)\rarr \Mod (S)$
such that, $\Delta(\lambda)^\circ \cong \nabla(\lambda)$ 
(and $\Ext^i(M,N) \cong \Ext^i(N^\circ ,M^\circ)$). 
We will usually shorten this and say $S$ has a simple preserving
duality.

\begin{rem}It is clear (given the equivalences in the definition for the
good filtration dimension) that for $S$ with simple preserving duality
and $M$ an $S$-module 
we have $\wfd(M) = \gfd(M^\circ)$.
We will use this without further comment.
\end{rem}


Schur algebras are quasi-hereditary with poset $\partn(n,r)$ ordered
by dominance~\cite[2.2h]{donk1}. 
The costandard modules correspond to the $\GL_n$-modules
$\nabla(\lambda)$.
Schur algebras  also have a simple preserving duality see~\cite[remark
(ii) following Lemma 4.1.3]{donkbk}. 
The rest of this paper is devoted to calculating
$\gfd\bigl(S_q(2,r)\bigr)$, $\glob\bigl(S_q(2,r)\bigr)$,
$\gfd\bigl(S(3,r)\bigr)$  and $\glob\bigl(S(3,r)\bigr)$
explicitly.

In all these cases we need the following technical lemma.
\begin{lem}\label{lem:gandi}
Suppose we have a short exact sequence of $S$-modules
$$0\rarr A\rarr B \rarr C \rarr 0. $$
\begin{enumerate}
\item[(i)]
\ If  $\gfd(B) > \gfd(C)$ then
$\gfd(A)=\gfd(B)$.
Furthermore if $\gfd(B) > \gfd(C)+1$ then for
all $S$-modules $M$ we have
$$\Ext^{\wfd(M)+\gfd(A)}(M, A) \cong \Ext^{\wfd(M)+\gfd(B)}(M, B).$$
\item[(ii)]
\ If  $\gfd(B) < \gfd(C)$ then
$\gfd(A)=\gfd(C)+1$.
Furthermore 
for all $S$-modules $M$ we have
$$\Ext^{\wfd(M)+\gfd(A)}(M, A) \cong \Ext^{\wfd(M)+\gfd(C)}(M, C).$$ 
\end{enumerate}
\end{lem}
\begin{proof}
We consider (i), the proof of (ii) is similar.
The long exact sequence corresponding to the short exact sequence gives us
$$\gfd(B) \le \max\{\gfd(A),\gfd(C)\}\mbox{ and } 
\gfd(A)\le \max\{\gfd(B),\gfd(C)+1\}.$$
So if $\gfd(B) > \gfd(C)$ then $\gfd(B)=\gfd(A)$.
If furthermore $\gfd(B) > \gfd(C)+1$ then for
all $S$-modules $M$ we have
$$\Ext^{\wfd(M)+\gfd(A)}(M, A) \cong \Ext^{\wfd(M)+\gfd(B)}(M, B)$$
using Lemma~\ref{lem:wandg} and the long exact sequence.
\end{proof}

\begin{lem}\label{lem:globtwo}
Suppose $S$ is a quasi-hereditary algebra with simple preserving
duality $^\circ$ and $\lambda \in \Lambda$  with $\gfd\bigl(L(\lambda)\bigr)=1$.
Then
$\Ext^2\bigl(L(\lambda),L(\lambda)\bigr)\cong \Hom(Q^\circ,Q)\not \cong
0$, where $Q=\nabla(\lambda)/L(\lambda)$. 
\end{lem}
\begin{proof}
We have a short exact sequence
$$0 \rarr Q^\circ \rarr \Delta(\lambda) \rarr L(\lambda) \rarr 0.$$
Dimension shifting gives us $\Ext^2\bigl(L(\lambda),
L(\lambda)\bigr)\cong \Ext^1\bigl(Q^\circ, L(\lambda)\bigr) $.
Now $\Hom\bigl(Q^\circ ,L(\lambda)\bigr)\cong
\Hom\bigl(Q^\circ,\break
\nabla(\lambda)\bigr) \cong0$ as $L(\lambda)$ is not a composition
factor of $Q$. 
So another dimension shift gives us 
 $\Ext^1\bigl(Q^\circ, L(\lambda)\bigr)\cong \Hom(Q^\circ, Q\bigr)$.
We can map
simple modules in the head of $Q^\circ$ to the socle of $Q$,
so $\Hom(Q^\circ, Q\bigr)\not \cong 0$ and we are done.
\end{proof}

\section{The $\SL_2$ case}
In this section we focus on $\SL_2$.
We will identify a block of a Schur algebra with the highest weights of the 
irreducibles that belong to the block. 

Suppose $\lambda \in X^+$ and write $\lambda=p \lambda_1
+\lambda_0$ with $\lambda_1\in X^+$ and $\lambda_0\in X_1$. 
We define $\nabla_p(\lambda)=\nabla(\lambda_1)\frob\otimes
L(\lambda_0)$, 
and $\Delta_p(\lambda)=\Delta(\lambda_1)\frob\otimes L(\lambda_0)$. 
It is clear that 
$\nabla_p(\lambda)=\bigl(\Delta_p(\lambda)\bigr)^\circ$. These modules
will play an important role in what follows.
We define $\hw(M)$ to be the set of weights $\lambda \in X^+$ with 
$L(\lambda)$ a composition factor of $M$ such that there is no $\mu \in X^+$ with
$\mu > \lambda$ and $L(\mu)$ a composition factor of $M$.

We now note a sequence,  
first remarked upon by Jantzen in~\cite[remark 2 following theorem
3.8]{jandar}
and proved by Xanthopoulos in~\cite[proposition 6.1.1]{xanth} 
which we will use repeatedly in this section.

\begin{lem}\label{lem:xanth}
We have for $r \ge 1$ and $0\le a \le p-2$ a
short exact sequence
$$0\rarr \nabla(a) \otimes \nabla(r)\frob \rarr\nabla(pr +a )
\rarr\nabla(p-a-2)\otimes \nabla(r-1) \frob \rarr 0.$$
If $a=p-1$ then $\nabla(a) \otimes \nabla(r)\frob \cong \nabla(pr +
a)$.
\end{lem}

We write $r=r_0 + pr_1$ where $0\le r_0 < p$.
If $\lambda=(r)$ then we define 
$g(\lambda)=r_1$.

We denote the block of $S(n,r)$ containing $(a_1,a_2,\ldots,a_n) \in\partn(n,r)$
by $\B(a_1,a_2,\ldots,a_n)$. For $n=2$, $\B(a_1,a_2)$ is totally ordered. 
A block $\B(a_1,a_2)$ is defined to be primitive if 
$a_1-a_2 \not\equiv -1 \pmod p$.
We say a weight $\lambda=(r)\in X^+$ is primitive if $r_0\ne p-1$. 

\begin{lem}\label{lem:Qweight}
Suppose $\lambda$ is primitive and  $Q$ is the quotient
$\nabla_p(\lambda)/L(\lambda)$.
Then $g\bigl(\hw(Q)\bigr)\le g(\lambda)-2$ with $g(\lambda)$ as
defined above.
\end{lem}
\begin{proof}
We write $\lambda=r = p r_1 + r_0$ and $r_1=pr_1^\prime +r_0^\prime$.
with $0\le r_0< p$ and $0\le r_0^\prime < p$.
Then $L(\lambda) \cong L(r_1)\frob \otimes L(r_0)$
by Steinberg's tensor product theorem and so $Q\cong (Q^\prime)\frob \otimes
L(r_0)$ where $Q^\prime$ is the quotient $\nabla(r_1)/L(r_1)$.
It is clear using~\ref{lem:xanth} that 
$[\nabla(r_1):L(pr_1^\prime -r_0^\prime-2)] \ne 0$. 
We have $[Q^\prime :L(r_1)]=0$.
Also the first weight in the $G$-block of $(r_1)$ smaller than $(r_1)$
is $(pr_1^\prime -r_0^\prime-2)$ using~\cite[main theorem]{donk4}.
Hence $\hw(Q^\prime) = (pr_1^\prime -r_0^\prime -2)$. Thus
$\hw(Q)=p(pr_1^\prime -r_0^\prime -2)+ (r_0)= (p(r_1-2r_0^\prime-2) + r_0)$.
Thus $g(\hw(Q))=r_1-2r_0^\prime-2\le g(\lambda)-2$ as required. 
\end{proof}

The next Theorem shows that $\gfd\bigl(L(\lambda)\bigr)=g(\lambda)$.
The following Lemma forms part of the inductive step.

\begin{lem}\label{lem:ind}
Let $\lambda \in X^+$ and 
suppose $g(\mu)=\gfd\bigl(L(\mu)\bigr)= \gfd\bigl(\nabla_p(\mu)\bigr)$ for
all $\mu < \lambda$ and $\gfd\bigl(\nabla_p(\lambda)\bigr)=g(\lambda)$. Then
$\gfd\bigl(L(\lambda)\bigr)=\gfd\bigl(\nabla_p(\lambda)\bigr)$.
Furthermore 
$$\Ext^{2 \gfd\left(L(\lambda)\right)}\bigl(L(\lambda),
L(\lambda)\bigr) \cong \Ext^{2
\gfd\left(L(\lambda)\right)}\bigl(\Delta_p(\lambda),\nabla_p(\lambda)\bigr).$$
\end{lem}
\begin{proof}
We have a short exact sequence
$$0\rarr L(\lambda) \rarr \nabla_p(\lambda) \rarr Q \rarr 0.$$
Lemma~\ref{lem:Qweight} gives us $\gfd\bigl(\nabla_p(\lambda)\bigr)> \gfd(Q)
+ 1$, as $\gfd(Q)$ is at most the maximum of the good filtration dimensions of its
composition factors.
The Lemma now follows by Lemma~\ref{lem:gandi} part (i) (applied
twice).
\end{proof}

\begin{thm}\label{thm:sl2}
Suppose $r=r_0+p r_1$ with $0\le r_0\le p-2$. Then
$$\gfd\bigl(\nabla(r_1)\frob\otimes L(r_0)\bigr)=\gfd\bigl(L(r)\bigr)=r_1$$
and
$$\Ext^{2 r_1}\bigl(L(r), L(r)\bigr)\cong
\Ext^{2 r_1}\bigl(\Delta_p(r), \nabla_p(r)\bigr)\cong k.$$ 
\end{thm}
\begin{proof} 
We proceed by induction on $r_1$.
For $r_1=0$ we have $L(r_0)=\nabla(r_0)$ and so
$\gfd\bigl(L(r_0)\bigr)=0$.
For $r_1=1$ we have $L(r)=\nabla(r_1)\frob\otimes\nabla(r_0)$
and a non-split short exact sequence 
$$0\rarr L(r)\rarr \nabla(r) \rarr \nabla(p-r_0-2)\rarr 0.$$
Thus $\gfd\bigl(L(r)\bigr)= 1$ by
Lemma~\ref{lem:globtwo}. This Lemma also gives us
$$\Ext^2\bigl(L(r),L(r)\bigr)\cong \Hom \bigl(L(p-r_0-2),
L(p-r_0-2)\bigr)\cong k.$$

We now suppose that $r_1\ge 2$.
Since by induction
$\gfd\bigl(\nabla(r_1-1)\frob\otimes L(p-r_0-2)\bigr)=r_1-1 \ge 1$, 
applying Lemma~\ref{lem:gandi} to the short exact sequence from
Lemma~\ref{lem:xanth} we have
$\gfd\bigl(\nabla(r_1)\frob\otimes L(r_0)\bigr)=r_1$.
By Lemma~\ref{lem:ind} and induction we have 
$\gfd\bigl(L(r)\bigr)= \gfd\bigl(\nabla_p(r)\bigr)=r_1$.

We also have $\Ext^{2r_1}\bigl(\Delta_p(r), \nabla_p(r)\bigr) \cong 
\Ext^{2r_1-2}\bigl(\Delta_p(p r_1 - r_0 -2), \nabla_p(p r_1 -r_0
-2)\bigr)$ by Lemma~\ref{lem:gandi} part (ii) and Lemma~\ref{lem:xanth}. 
But by induction this is isomorphic to $k$.
Lemma~\ref{lem:ind} then completes the proof.
\end{proof}

\begin{cor}\label{cor:sl2}
Suppose $r$ is not primitive and write $r=p^d r_1+ p^d -1$ for some $d \in
\N^+$ with $r_1\not\equiv-1 \pmod p$. 
Then $\gfd\bigl(L(r)\bigr)=\gfd\bigl(L(r_1)\bigr)$ and 
$$\Ext^{2 \gfd(L(r))}\bigl(L(r), L(r)\bigr)\cong k.$$
\end{cor}
\begin{proof}
We have by~\cite[section 4, theorem]{donk4} that 
$\B(r,0)$ is Morita equivalent to $\B(r_1,0)$ in $S(2,r_1)$ with $r_1$
primitive and the result follows.
\end{proof}

\begin{cor}\label{cor:S2gfd}
Given $(a_1,a_2)\in \partn(2,r)$, we write $a_1-a_2= p^{d+1}c_1+p^dc_0+p^d-1$
with $d\in\N$ and $0\le c_0\le p-2$.
Then $\gfd\bigl(L(a_1,a_2)\bigr)=c_1$ where $L(a_1,a_2)$ is the irreducible
module of highest weight $(a_1,a_2)$ for $S(2,r)$.
Moreover
$$\Ext^{2c_1}_{S(2,r)}\bigl(L(a_1,a_2),L(a_1,a_2)\bigr)\cong k.$$
\end{cor}
\begin{proof}
Now Theorem~\ref{thm:sl2} and Corollary~\ref{cor:sl2} give us
$$\Ext_{S(2,r)}^i\bigl(\Delta(b_1,b_2), L(a_1,a_2)\bigr)
\cong\Ext_{G}^i\bigl(\Delta(b_1-b_2), L(a_1-a_2)\bigr)
\cong 0 
$$
if $i> c_1$
and so we have $\gfd\bigl(L(a_1,a_2)\bigr) \le c_1$.
We also have 
$$\Ext^{2c_1}_{S(2,r)}\bigl(L(a_1,a_2),L(a_1,a_2)\bigr)
\cong\Ext^{2c_1}_G\bigl(L(a_1-a_2),L(a_1-a_2)\bigr)
\cong k.$$
Thus Lemma~\ref{lem:wandg} gives us
$\wfd L(a_1,a_2)+\gfd L(a_1,a_2)\ge 2 c_1$,
But $\wfd L(a_1,a_2)=\gfd L(a_1,a_2)$, and 
so we have $\gfd L(a_1,a_2)= c_1$, 
as required.
\end{proof}

\begin{thm}\label{thm:s2glob}
The global dimension of $S(2,r)$ is twice its good filtration
dimension and is given as follows:
$$\begin{array}{llrcl}
\mbox{for}& p=2 & \glob(S(2,r)) &= &\left\{ \begin{array}{ll}
                                r &\mbox{if $r$ is even}\cr
                                2\lfloor\frac{r}{4}\rfloor & \mbox{if
                                        $r$ is odd}
                                \end{array}
                        \right.\\
&  & & & \\
\mbox{for} & p\ge3 & \glob(S(2,r)) &= &2\lfloor\frac{r}{p}\rfloor.
\end{array}$$
\end{thm}
\begin{proof}
If $r$ is even then $\Lambda^+(2,r)$ consists of all partitions
$(a,b)$ of $r$ whose difference $a-b$ are all even numbers less
than~$r$. 
If $r$ is odd then $\Lambda^+(2,r)$ consists of all partitions
$(a,b)$ of $r$ whose difference $a-b$ are all odd numbers less 
than~$r$. 

We let $S=S(2,r)$.
Using Corollary~\ref{cor:glob} and simple preserving duality we know 
$\glob(S) \le 2 \gfd(S)$. 
We also have 
$$\gfd(S) = \max\{ \gfd\bigl(L(a_1,a_2)\bigr) \mid
(a_1,a_2) \in \partn(2,r) \}.$$
But Corollary~\ref{cor:S2gfd} gives us 
$$\glob(S) \ge 2 \max\{ \gfd\bigl(L(a_1,a_2)\bigr) \mid
(a_1,a_2) \in \partn(2,r) \}.$$
Hence $\glob(S) = 2 \gfd(S)$.

Using the main theorem in~\cite{donk4} we have
for $p=2$ and $r$ even that $\B(r,0)=\partn(2,r)$ is primitive
so $\glob(S)=2 \frac{r}{2}=r$, as all other simple modules
corresponding to weights  in
$\partn(2,r)$ have smaller or equal good filtration dimension to that
of $L(r,0)$.
For $p=2$ and $r$ odd, $\B(r,0)$ is not primitive. We write $r= 2^d r_1+
2^d -1$ for some $d \in \N$ and $r_1$ even.
If $d=1$ then all other  simples modules corresponding to weights in
$\partn(2,r)$ have smaller or
equal good filtration dimension to that of $L(r,0)$, so $\glob(S)= r_1 =
2\lfloor\frac{r}{4}\rfloor$. 
If $d\ge 2$ then we consider $r^\prime=r-2$ so $(r^\prime+1,1)\in
\partn(2,r)$.
We have $r^\prime=2( 2^{d-1}(r_1-2))+1$ with $2^{d-1}(r_1-2)$ even,
  so $\glob(\B(r^\prime+1,1))=\glob(\B(r^\prime,0))=
 2^{d-1}(r_1-2)= \frac{r-3}{2}= 2\lfloor\frac{r}{4}\rfloor$.

For $p\ge3$ and $r$ primitive then
$\glob(S)=2r_1=2\lfloor\frac{r}{p}\rfloor$,
as all other simple modules corresponding to weights  in $\partn(2,r)$ 
have smaller or equal good filtration dimension to that of $L(r,0)$.
If $r$ is not primitive then $r= p^d r_1+
p^d -1$ for some $d \in \N$ and $r_1\not\equiv-1 \pmod p$. 
Consider $r^\prime = r-2$, then $(r^\prime+1,1) \in
\partn(2,r)$ is primitive and so $\glob(S)=2 \lfloor\frac{r^\prime}{p}\rfloor
=\lfloor\frac{r}{p}\rfloor$.
\end{proof}

We now consider the quantum version of the Schur algebra 
$S_q(n,r)$ with $0 \ne q\in k$. This is a deformation of
the classical Schur algebra with parameter $q$.
See the introduction
of~\cite{donkbk} for the basic properties of $S_q(n,r)$.
When $q=1$ then $S_q(n,r)$ is just the classical Schur algebra.
If $q$ is not a root of unity then $S_q(n,r)$ is semi-simple.
We will consider the case where $q$ is a primitive $l$th root of unity with 
$l \ge 2$.


We now show that the argument above generalises to the quantum case.
To do this we need the appropriate quantum versions of the results used
above. We will be using the Dipper-Donkin quantum group $q$-$\GL_n$ defined 
in~\cite{dipdonk}.  

Now we know that $S_q(n,r)$ is quasi-hereditary 
with poset $\partn(n,r)$ by~\cite[section 4, (6)]{donkquant}.
We also have the property for all $S_q(n,r)$-modules $V$ and $W$
that
$$\Ext^i_{S_q(n,r)}(V,W) \cong \Ext^i_{q\mbox-\GL_n}(V,W) $$
by~\cite[section 4, (5)]{donkquant}.
The blocks of $S_q(n,r)$ were 
determined in~\cite{cox2}. We also know that all blocks of
$S_q(n,r)$ are Morita equivalent to a block of $S_q(n,r^\prime)$ with
$r^\prime$ primitive by~\cite[lemma 6.10]{coxerd}.
We have a quantum Frobenius morphism $\mathrm{F}: q\mbox-\GL_n \rarr
GL_n(k)$.
Some of the other basic properties of $q$-$\GL_n$-modules appears 
in~\cite[chapter 3]{donkbk} including a proof of the quantum version of
Steinberg's tensor product theorem.
Suppose we write $\lambda= l\lambda_1 + \lambda_0$ with $\lambda_0$
$l$-restricted and $\lambda_1$ dominant. We define $\nabla_l(\lambda)=
\overline{\nabla}(\lambda_1)\frob \otimes L(\lambda_0)$, where
$\overline{\nabla}(\lambda_1)$ is the classical module in characteristic $p$.

We now let $n=2$. 
The generalisation of Lemma~\ref{lem:xanth} appears
in~\cite[Proposition 3.1]{cox1}.
Moreover all it does is to replace
$p$ with $l$. Thus Lemmas~\ref{lem:Qweight} and~\ref{lem:ind}
will carry through unchanged.
Theorem~\ref{thm:sl2} and Corollary~\ref{cor:S2gfd} now generalise to give:

\begin{thm}\label{thm:qgl2}
Suppose $r=r_0+l r_1$ with $0\le r_0\le l-2$. Then
$$\gfd\bigl(\overline{\nabla}(r_1,0)\frob\otimes L(r_0,0)\bigr)
=\gfd\bigl(L(r,0)\bigr)=r_1.$$
and
$$\Ext^{2 r_1}\bigl(L(r,0), L(r,0)\bigr)\cong
\Ext^{2 r_1}\bigl(\Delta_l(r,0), \nabla_l(r,0)\bigr)\cong k.$$ 
\end{thm}

Also we can use~\cite[lemma 6.10]{coxerd} and the classical
result~\ref{thm:sl2} to give the good filtration
dimension of $L(\lambda)$ for $\lambda$ non-primitive, so
we have the following Corollary.

\begin{cor}\label{cor:qsl2}
Suppose $r$ is not primitive so $r=lp^{d-1} r_1+ lp^{d-1} -1$ for some $d \in
\N^+$ and  $r_1\not\equiv-1 \pmod p$. 
Then
$\gfd\bigl(L(r,0)\bigr)=\gfd\bigl(\bar{L}(r_1,0)\bigr)
=\lfloor\frac{r_1}{p}\rfloor$ and 
$$\Ext^{2 \gfd(L(r,0))}\bigl(L(r,0), L(r,0)\bigr)\cong k$$
where $\bar{L}(r_1,0)$ is the irreducible $S(2,r_1)$-module of highest weight 
$(r_1,0)$.
\end{cor}

We now have the following theorem for the quantum case.

\begin{thm}
The global dimension of $S_q(2,r)$ is twice its good filtration
dimension and is given as follows:
$$\begin{array}{llrcl}
\mbox{for}& l=2 & \glob(S_q(2,r)) &= &\left\{ \begin{array}{ll}
                                r &\mbox{if $r$ is even}\cr
                                2\lfloor\frac{r}{2p}\rfloor & \mbox{if
                                        $r$ is odd}
                                \end{array}
                        \right.\\
&  & & & \\
\mbox{for} & l\ge3 & \glob(S_q(2,r)) &= &2\lfloor\frac{r}{l}\rfloor.
\end{array}$$
\end{thm}
\begin{proof}
%
The argument is very similar to that of Theorem~\ref{thm:s2glob}.
%
%
%
%
\end{proof}

\section{Some filtrations for $\SL_3$}
In this section we obtain a filtration of the modules
$\nabla(\lambda)$ for $\SL_3$  by modules of
the form $\nabla(\lambda)\frob \otimes L(\mu)$ with $\lambda \in X^+$
and $\mu \in X_1$. We call such a filtration a \emph{$p$-filtration}.
Throughout this section $G=\SL_3$.
We start by proving some Lemmas about extensions between two modules 
$\nabla_p(\lambda)$ and $\nabla_p(\mu)$.

Given a rational $G$-module $V$, we have a five term exact sequence,
\begin{eqnarray*}
\lefteqn{
0\rarr H^1(G/{G_1}, V^{G_1})\rarr H^1(G,V)\rarr H^1(G_1,V)^{G/G_1} 
} & \hspace{4cm} & \\
& & 
\rarr H^2(G/{G_1},V^{G_1}) \rarr H^2(G,V).
\end{eqnarray*}
\hspace{-8pt} 
This is the Lyndon-Hochschild-Serre sequence for $G$ and $G_1$.
A $G$-module $W$ which is trivial as a $G_1$-module, is of the form $V\frob$ for some
$G$-module $V$ which is unique up to isomorphism. We define $W^{(-1)}=V$.
If $W$ and $V$ are $G$-modules then $\Ext^i_{G_1}(W,V)$ has a natural
structure as a $G$-module.
Moreover when $W$ and $V$ are finite dimensional we have,
\begin{eqnarray}
\label{frob}
\Ext^i_{G_1}(W,V\otimes Y\frob)\cong \Ext^i_{G_1}(W,V) \otimes Y\frob
\end{eqnarray}
\hspace{-7pt}
as $G$-modules.
If $H=G$ or $G_1$ then we have
\begin{eqnarray}
\label{extH}
\Ext^i_H(W,V)\cong \Ext^i_H(V^*,W^*)\cong \Ext^i_H(k,W^*\otimes V) \cong
H^i(H,W^*\otimes V)
\end{eqnarray}
\hspace{-7pt}
where $^*$ is the ordinary dual. We have $\nabla(a,b)^*\cong
\Delta(b,a)$.
We also note that $(V\frob)^{G/G_1}\cong
V^G$, and $H^i(G/G_1,V\frob) \cong H^i(G,V)$.

In the following sections we make repeated use of a Proposition proved
in the PhD thesis of Yehia. We reproduce his results here for the 
convenience of the reader.
\begin{propn}\label{propn:yehia}
The non-zero $\Ext^1_{G_1}\bigl(L(\alpha),L(\beta)\bigr)$ for $\alpha$, $\beta \in
X_1$ are given by the following tables.
\begin{enumerate}
\item[(i)]\ 
For $(r,s) \in X_1$ with $r+s=p-2$, we have
\begin{center}
\begin{tabular}{c|ccc} 
$\alpha \downarrow$, $\beta \rightarrow$&$(r,s)$& $(p-1,r)$ & $(s,p-1)$ \\\hline 
$(r,s)$   &$0$    &$\nabla(0,1)\frob$&$\nabla(1,0)\frob$ \\
$(p-1,r)$ &$\nabla(1,0)\frob$& $0$     & $0$ \\
$(s,p-1)$ &$\nabla(0,1)\frob$& $0$     & $0$ \\
\end{tabular}
\end{center}
\medskip
\item[(ii)]\ 
For $(r,s) \in A_0$ and $ p \ge 5$, we have
\begin{center}
\begin{tabular}{c|ccc}
$\alpha \downarrow$, $\beta \rightarrow$ &$(p-s-2,p-r-2)$   
& $(r+s+1,p-s-2)$ & $(p-r-2,r+s+1)$ \\\hline 
$(r,s)$   &$k$    &$\nabla(0,1)\frob$&$\nabla(1,0)\frob$ \\
\end{tabular}
\end{center}
\medskip
\begin{center}
\begin{tabular}{c|ccc} 
$\alpha \downarrow$, $\beta \rightarrow$&$(r,s)$   
& $(s,p-r-s-3)$ & $(p-r-s-3,r)$ \\\hline 
$(p-s-2,p-r-2)$   &$k$               
&$\nabla(0,1)\frob$&$\nabla(1,0)\frob$ \\
\end{tabular}
\end{center}
\medskip
If $p=3$ then all the entries in the two tables above are replaced by
$k \oplus \nabla(0,1)\frob\oplus \nabla(1,0)\frob$.
\end{enumerate}
\end{propn}
\begin{proof}
See~\cite[proposition 3.3.2]{yehia}.
\end{proof}

\begin{lem}\label{lem:extone}
For $\lambda$, $\mu \in X^+$ and $\alpha$, $\beta \in X_1$ we have,
$$\begin{array}{l}
\Ext^1_G(\nabla(\lambda)\frob \otimes L(\alpha),
\nabla(\mu)\frob\otimes L(\beta)) \\ 
\hspace{15pt} \cong \left\{ \begin{array}{ll}
                        \Ext^1_G\bigl(\nabla(\lambda),\nabla(\mu)\bigr) 
                                &\mbox{if $\alpha=\beta$,}\\
                        \Hom_G\left(\nabla(\lambda),\nabla(\mu)\otimes
                        \Ext^1_{G_1}\bigl(L(\alpha),L(\beta)\bigr)^{(-1)}\right)
                                  &\rlap{\mbox{otherwise.}}  
                        \end{array}\right.
\end{array}$$
\end{lem}
\begin{proof}
Case (i) $\alpha=\beta$. We know by Proposition~\ref{propn:yehia} that 
$\Ext_{G_1}^1(L(\alpha), L(\beta))=0$, so by~(\ref{frob}) and~(\ref{extH})
we have
$$\left(H^1(G_1,\Delta(\lambda^*)\frob\otimes\nabla(\mu)\frob\otimes
L(\beta)\otimes L(\alpha^*))\right)^{G/{G_1}} = 0.$$
Thus~(\ref{extH}) and the five term exact sequence gives us:
$$\begin{array}{rl}
\Ext^1_G(\nabla(\lambda)\frob \otimes L(\alpha),
\nabla(\mu)\frob\otimes L(\beta))\hspace{-80pt} \\ 
&\cong H^1(G,\Delta(\lambda^*)\frob\otimes\nabla(\mu)\frob\otimes
L(\beta)\otimes L(\alpha^*))\\
&\cong H^1\Bigl(G/{G_1}, \left(\Delta(\lambda^*)\frob
\otimes\nabla(\mu)\frob\otimes L(\beta)\otimes
L(\alpha^*)\right)^{G_1}\Bigr).
\end{array}$$
But $\alpha=\beta$ so  
$$\left(\Delta(\lambda^*)\frob\otimes\nabla(\mu)\frob\otimes
L(\beta)\otimes L(\alpha^*)\right)^{G_1} = \Delta(\lambda^*)\frob
\otimes\nabla(\mu)\frob.$$
and hence
$$\begin{array}{rl}
\Ext^1_G(\nabla(\lambda)\frob \otimes L(\alpha),
\nabla(\mu)\frob\otimes L(\beta))  
& \cong H^1(G/{G_1}, \Delta(\lambda^*)\frob
        \otimes\nabla(\mu)\frob)\\
& \cong H^1(G, \Delta(\lambda^*)
        \otimes\nabla(\mu))\\
& \cong \Ext^1_G(\nabla(\lambda),\nabla(\mu)\bigr)
\end{array}$$
as required.

Case (ii) $\alpha\ne \beta$.
We have 
$$\begin{array}{rl}
\Ext^1_G(\nabla(\lambda)\frob \otimes L(\alpha),
\nabla(\mu)\frob\otimes L(\beta))  
\hspace{-80pt}\\ & 
\cong H^1(G,\Delta(\lambda^*)\frob\otimes\nabla(\mu)\frob\otimes
L(\beta)\otimes L(\alpha^*)).
\end{array}
$$
In this case 
$$\left(\Delta(\lambda^*)\frob\otimes\nabla(\mu)\frob\otimes
L(\beta)\otimes L(\alpha^*)\right)^{G_1} = 0$$
as $L(\beta)\otimes L(\alpha^*)$ has no $G_1$ fixed points.
So by the five term exact sequence we have
$$\begin{array}{rl}
\Ext^1_G(\nabla(\lambda)\frob \otimes L(\alpha),
\nabla(\mu)\frob\otimes L(\beta))\hspace{-140pt} \\ 
& \cong \left(H^1(G_1,\Delta(\lambda^*)\frob\otimes\nabla(\mu)\frob\otimes
L(\beta)\otimes L(\alpha^*))\right)^{G/{G_1}}\\
& \cong
        \left(\Delta(\lambda^*)\frob\otimes\nabla(\mu)\frob\otimes
        \Ext^1_{G_1}(L(\alpha),L(\beta))\right)^{G/{G_1}}\\
&\cong 
        \left(\Delta(\lambda^*)\otimes\nabla(\mu)\otimes
        \Ext^1_{G_1}(L(\alpha),L(\beta))^{(-1)}\right)^{G}\\
&\cong
        \Hom_G\Bigl(\nabla(\lambda),\nabla(\mu)\otimes
\left(\Ext^1_{G_1}\bigl(L(\alpha),L(\beta)\bigr)\right)^{(-1)}\Bigr)
\end{array}$$
as required. \end{proof}

\begin{lem}\label{lem:exta}
Let $(a-1,b+1),(a,b-1) \in X^+$ and $(r,s)\in X_1$, then
$$\begin{array}{l}
\Ext^1_G\bigl(\nabla(a-1,b+1)\frob\otimes L(r,s), 
\nabla(a,b-1)\frob\otimes L(r,s)\bigr) \\
\hspace{150pt} \cong\left\{\begin{array}{ll}
      0 &\mbox{if $b \not\equiv -1 \pmod p$}\\
      k &\mbox{if $b \equiv -1 \pmod p$}.
    \end{array}\right. 
\end{array}$$
\end{lem}
\begin{proof}
Lemma~\ref{lem:extone} gives us 
$$\begin{array}{l}
\Ext^1_G\bigl(\nabla(a-1,b+1)\frob\otimes L(r,s), 
\nabla(a,b-1)\frob\otimes L(r,s)\bigr) \\
\cong 
\Ext^1\bigl(\nabla(a-1,b+1),\nabla(a,b-1)\bigr).\\
\end{array}$$
Now $(a-1,b+1)$ and  $(a,b-1)$ differ by a
single root so we may apply the result of~\cite[(4.3) and (3.6)]{erd1} and we
are done.
\end{proof}

Similarly we have:

\begin{lem}\label{lem:extb}
Let $(a+1,b-1),(a-1,b) \in X^+$ and $(r,s)\in X_1$, then
$$\Ext^1_G\bigl(\nabla(a+1,b-1)\frob\otimes L(r,s),
\nabla(a-1,b)\frob\otimes L(r,s)\bigr)\cong\left\{\begin{array}{ll}
                                0 &\mbox{if $a \not\equiv -1 \pmod p$}\\
                                k &\mbox{if $a \equiv -1 \pmod p$}.
                                        \end{array}\right.
$$
\end{lem}

We also have:
\begin{lem}\label{lem:extzero}
Suppose $\alpha, \beta \in X_1$ with $\alpha\ne \beta$,
and $\alpha$ and $\beta$ in the same alcove.
Then for all $\lambda, \mu \in X^+$ we have 
$$\Ext^1\bigl(\nabla(\lambda)\frob\otimes L(\alpha), \nabla(\mu)\frob\otimes
L(\beta)\bigr) \cong 0.$$
\end{lem}
\begin{proof}
By Proposition~\ref{propn:yehia} we have
$$\Ext^1_{G_1}\bigl(L(\alpha),L(\beta)\bigr)\cong 0$$
and the result follows by Lemma~\ref{lem:extone}.
\end{proof}

In the main theorem in this section we use repeatedly the fact that
the modules $\nabla(\lambda)$ have both simple socle and simple head
for $\lambda$ dominant. Jantzen proved that $\nabla(\lambda)$ has
simple head for $\SL_3$ when $p >3$. The following four Lemmas extend
this result to $p=2$ and $3$. In these Lemmas $\rho =(1,1)$,
which is half the sum of all positive roots for the Weyl group, 
$w_0$ is the longest word in the Weyl group
and $\St=L(p-1,p-1)$ is the Steinberg module.
We also have a $G$-module 
$Q_1(\mu)$ which when considered as a $\uone$-$T$-module is the 
injective hull of $L(\mu)$ for $\mu\in X_1$, 

\begin{lem}\label{lem:good}
For all $\lambda$ and $\mu \in X^+$ the module
$\nabla(\lambda) \otimes \nabla(\mu)$ has a good filtration.
Moreover the $\nabla(\nu)$ which appear as quotients in this
filtration are given by Brauer's character formula.
\end{lem}
\begin{proof}
A proof of the property that $\nabla(\lambda) \otimes \nabla(\mu)$ has 
a good filtration, for type $A_n$, is given in~\cite{wang1}. 
It is proved for most other cases in~\cite{donkrat}. The general proof
is given in~\cite{mathieu}.
For a version of Brauer's character formula see~\cite[II, lemma
5.8]{jantz}.
\end{proof}

\begin{cor}
$\St$ is a direct summand of $\St \otimes \St$.
\end{cor}
\begin{proof}
All the quotients $\nabla(\nu)=\nabla(\nu_1,\nu_2)$ which appear in 
a good filtration of
$\St\otimes \St$ have $0 \le \nu_i \le 2p-2)$. Also $(0,0)$ is a weight
in $\St$. Hence using Brauer's character formula, $\St=\nabla(p-1,p-1)$
appears as a quotient in a good filtration of $\St\otimes \St$, but
all other weights that appear are not linked to $(p-1,p-1)$. (The
first such possible weight is $(2p-1,2p-1)$.) Thus
$\St$ is a direct summand of $\St\otimes \St$ as required.
\end{proof}

\begin{lem}\label{lem:Qweyl}
If $\lambda\in X^+$ and $\mu\in X_1$ then $\Delta(\lambda)\frob
\otimes Q_1\bigl(p\rho+w_0(\mu+\rho)\bigr)$ has a Weyl filtration.
\end{lem}
\begin{proof}
This is proved for $p>3$ in~\cite[theorem 5.6]{jandar}.
Now~\cite[lemma 3.1.3]{yehia} shows for $p=2$ and $3$  
that $Q_1\bigl(p\rho+w_0(\mu+\rho)\bigr)$ 
is a $G$-direct summand of $L(\mu)\otimes \St$.
Hence $\Delta(\lambda)\frob
\otimes Q_1\bigl(p\rho+w_0(\mu+\rho)\bigr)$ is a direct summand of
$\Delta(\lambda)\frob \otimes L(\mu)\otimes \St$ which in turn is a
direct summand of $\Delta(\lambda)\frob \otimes L(\mu)\otimes
\St\otimes \St$. But this last module has a Weyl filtration. First
$\Delta(\lambda)\frob \otimes \St \cong
\Delta\bigl(p\lambda+(p-1,p-1)\bigr)$,
also $L(\mu)\otimes \St$ has a Weyl filtration by~\cite[lemma 3.1.2]{yehia}.
Hence $\Delta(\lambda)\frob \otimes L(\mu)\otimes
\St\otimes \St$ has a Weyl filtration and we are done.
\end{proof}

\begin{cor}\label{cor:Qweyl}
Suppose $\lambda \in X^+$ and $\mu\in X_1$ then
$\Delta\bigl(p(\lambda+\rho)-\rho+\mu\bigr)$ embeds in
$\Delta(\lambda)\frob\otimes Q_1\bigl(p\rho+w_0(\mu+\rho)\bigr)$.
\end{cor}
\begin{proof}
This follows using characters as in the proof of~\cite[corollary
5.7]{jandar}.
\end{proof}

\begin{lem}\label{lem:weylsoc}
For all $\lambda\in X^+$ and $\mu\in X_1$ we have
$$\soc_{G}\Delta\bigl(p(\lambda+\rho)+\mu-\rho\bigr)
\cong \bigl(\soc_G\Delta(\lambda)\bigr)\frob\otimes 
L\bigl(p\rho+w_0(\mu+\rho)\bigr).$$
\end{lem}
\begin{proof}
This follows as in Jantzen~\cite[theorem 6.2]{jandar} using
Corollary~\ref{cor:Qweyl}
to remove the restriction on $p$.
\end{proof}

\begin{propn}\label{propn:simpsoc}
For all $\lambda\in X^+$ the module $\Delta(\lambda)$ has simple socle.
\end{propn}
\begin{proof}
We can use~\ref{lem:weylsoc} to remove the restriction on $p$ in 
the argument of Jantzen in~\cite[6.9]{jandar}, the result then
follows.
\end{proof}

We will now give explicitly the $p$-filtrations of $\nabla(\lambda)$. 
Further information about these filtrations (for $p > 3$) appears 
in~\cite[chapter 2, lemma 2 and following]{kuehne}.
In the following diagrams the sections of a filtration of
$\nabla(\lambda)$ are written as modules with
lines connecting them. Sections that have a non-trivial extension
between them that appear in $\nabla(\lambda)$ are joined by a line.
Sections that only have trivial extensions appearing between them are not
joined by a line.
Sections that embed in $\nabla(\lambda)$ are at the bottom of the
diagram. The sections that occur above it, are above these sections,
and so on.
In what follows if $\lambda_1$ is not dominant but one of its parts is
$-1$ then we take $\nabla_p(\lambda)$ to be the zero module. In other
words, it does not appear as a section in a $p$-filtration.

\begin{thm}\label{thm:pfiltra}
Each $\nabla(\lambda)$ has a $p$-filtration.
This filtration takes the following form:
\begin{enumerate}
\item[(i)]\ Suppose 
$\lambda = p(a,b) + (p-1,p-1)$ with $(a,b) \in X^+$. Then
$$\nabla(\lambda)= \nabla(a,b)\frob\otimes L(p-1,p-1).$$
\item[(ii)]\ Suppose
$\lambda = p(a,b) + (p-1,r)$ with $(a,b) \in X^+$ and $(p-1,r)\in
X_1$.
If we set $s=p-r-2$ then for $a \equiv -1 \pmod p$,
the module $\nabla(\lambda)$ has filtration 
$$\xymatrix@R=10pt
{
{\nabla(a,b-1)\frob\otimes L(s,p-1)} \ar@{-}[d] \\
{\nabla(a+1,b-1)\frob\otimes L(r,s)} \ar@{-}[d] \\
{\nabla(a-1,b)\frob\otimes L(r,s)} \ar@{-}[d] \\
{\nabla(a,b)\frob\otimes L(p-1,r)}}
$$ 
while for  $a \not\equiv -1 \pmod p$,
$\nabla(\lambda)$ has filtration 
$$\xymatrix@R=10pt@C=-20pt
{
& {\nabla(a,b-1)\frob\otimes L(s,p-1)} \ar@{-}[dl] \ar@{-}[dr]
\\
{\nabla(a+1,b-1)\frob\otimes L(r,s)} \ar@{-}[dr] &
 &{\nabla(a-1,b)\frob\otimes L(r,s)} \ar@{-}[dl] \\
& {\nabla(a,b)\frob\otimes L(p-1,r).}}
$$ 
\item[(iii)]\  Suppose
$\lambda = p(a,b) + (s,p-1)$ with $(a,b) \in X^+$ and $(s,p-1)\in
X_1$.
If we set $r=p-s-2$ then for $b \equiv -1 \pmod p$, the module
$\nabla(\lambda)$ has filtration 
$$\xymatrix@R=10pt
{
{\nabla(a-1,b)\frob\otimes L(p-1,r)} \ar@{-}[d] \\
{\nabla(a-1,b+1)\frob\otimes L(r,s)} \ar@{-}[d] \\
{\nabla(a,b-1)\frob\otimes L(r,s)} \ar@{-}[d] \\
{\nabla(a,b)\frob\otimes L(s,p-1)}}
$$ 
while for  $b \not\equiv -1 \pmod p$,
$\nabla(\lambda)$ has filtration 
$$\xymatrix@R=10pt@C=-20pt
{
& {\nabla(a-1,b)\frob\otimes L(p-1,r)} \ar@{-}[dl] \ar@{-}[dr]
\\
{\nabla(a-1,b+1)\frob\otimes L(r,s)} \ar@{-}[dr] &
 &{\nabla(a,b-1)\frob\otimes L(r,s)} \ar@{-}[dl] \\
& {\nabla(a,b)\frob\otimes L(s,p-1).}}
$$ 
\item[(iv)]\ Suppose
$\lambda = p(a,b) + (r,s)$ with $(a,b) \in X^+$, $(r,s)\in X_1$ and
$r+s=p-2$. Then the module
$\nabla(\lambda)$ has filtration
$$\xymatrix@R=10pt@C=-20pt
{
& {\nabla(a-1,b-1)\frob\otimes L(r,s)} \ar@{-}[dl] \ar@{-}[dr]
\\
{\nabla(a,b-1)\frob\otimes L(p-1,r)} \ar@{-}[dr] &
 &{\nabla(a-1,b)\frob\otimes L(s,p-1)} \ar@{-}[dl] \\
& {\nabla(a,b)\frob\otimes L(r,s).}}
$$
\item[(v)]\ Suppose 
$\lambda = p(a,0) + (r,s)$ with $(a,0) \in X^+$, $a\ge1$ and $(r,s)\in A_0$ 
then the module $\nabla(\lambda)$ has filtration
$$\xymatrix@R=10pt{
{\nabla(a-2,0)\frob\otimes L(s,p-r-s-3)} \ar@{-}[d] \\
{\nabla(a-1,0)\frob\otimes L(p-r-2,r+s+1)} \ar@{-}[d] \\
{\nabla(a,0)\frob\otimes L(r,s).}}
$$
\item[(vi)]\ Suppose
$\lambda = p(0,b) + (r,s)$ with $(0,b) \in X^+$, $b\ge1$ and $(r,s)\in A_0$.
Then the module $\nabla(\lambda)$ has filtration
$$\xymatrix@R=10pt{
{\nabla(0,b-2)\frob\otimes L(p-r-s-3,r)} \ar@{-}[d] \\
{\nabla(0,b-1)\frob\otimes L(r+s+1,p-s-2)} \ar@{-}[d] \\
{\nabla(0,b)\frob\otimes L(r,s).}}
$$
\item[(vii)]\ Suppose
$\lambda = p(a,b) + (r,s)$ with $(a,b) \in X^+$, $a$ and $b\ge 1$, 
and $(r,s)\in A_0$. We let 
$$\begin{array}{rlrl}\mu_1&=\lambda, &\mu_2&=(p a+r+s+1,p b-s-2),\\
\mu_3&=(pa +p-r-s-3,p b-2 p+r), &\mu_4&=(p a-r-2,p b+r+s+1),\\
\mu_5&=(p a-2 p+s,p b +p-r-s-3), &\mu_6&=(p a+s,p b-r-s-3),\\
\mu_7&=(p a-p+r,p b-p+s), &\mu_8&=(pa -r-s-3,p b+r), \\ 
\mu_9&=(p a-s-2,p b-r-2).& 
\end{array}$$
These weights are depicted in Figure~1~(a), where the number corresponds to
the subscript of $\mu$.
\begin{figure}[ht]
\begin{center}
\epsfbox{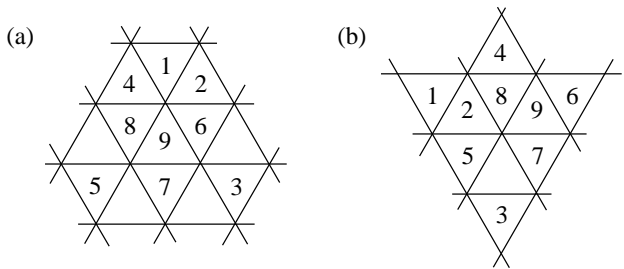}
\end{center}
\caption{Diagram showing weights for $\lambda$ inside (a) a lower
alcove and (b) an upper alcove}
\end{figure}

Then for $a$ and $b \equiv 0
\pmod p$, $\nabla(\lambda)$ has filtration 
$$\xymatrix@R=10pt@C=5pt
{
& &{\nabla_p(\mu_9)} \ar@{-}[dll] 
\ar@{-}[drr] \ar@{-}[d] \\
{{\nabla_p(\mu_6)}} \ar@{-}[ddrrr]   \ar@{-}[d]   
& & {\nabla_p(\mu_7)} \ar@{-}[ddr] \ar@{-}[ddl] 
& & {\nabla_p(\mu_8)} \ar@{-}[ddlll] \ar@{-}[d] \\
{\nabla_p(\mu_5)} \ar@{-}[dr]
& & & &{\nabla_p(\mu_3)} \ar@{-}[dl] \\
& {\nabla_p(\mu_4)} \ar@{-}[dr]
& &{\nabla_p(\mu_2)} \ar@{-}[dl] \\
& &{\nabla_p(\mu_1).}}
$$
For $a  \not\equiv 0 \pmod p$ there is no extension of
$\nabla_p(\mu_5)$ by $\nabla_p(\mu_6)$.
For $b  \not\equiv 0 \pmod p$ there is no extension of
$\nabla_p(\mu_3)$ by $\nabla_p(\mu_8)$.
So for $a$ and $b \not\equiv 0 \pmod p$ we have:
$$\xymatrix@R=10pt@C=5pt
{
& & &{\nabla_p(\mu_9)} \ar@{-}[dlll] \ar@{-}[drrr] \ar@{-}[dll] 
                \ar@{-}[drr] \ar@{-}[d] \\
{\nabla_p(\mu_5)} \ar@{-}[drr]
& {\ \nabla_p(\mu_6)\!\!\!} \ar@{-}[drrr] \ar@{-}[dr]    
& & {\nabla_p(\mu_7)} \ar@{-}[dr] \ar@{-}[dl] 
& & {\!\!\!\nabla_p(\mu_8)\ } \ar@{-}[dlll] \ar@{-}[dl] 
& {\nabla_p(\mu_3)} \ar@{-}[dll] \\
& &{{\!\nabla_p(\mu_4)\!}} \ar@{-}[dr]
& & {{\!\nabla_p(\mu_2)\!}} \ar@{-}[dl] \\
& & &{\nabla_p(\mu_1)}}
$$
and similarly for the other cases for $a$ and $b$.
\item[(viii)]\ Suppose
$\lambda = p(a,b) + (p-s-2,p-r-2)$ with $(a,b) \in X^+$, 
and $(r,s)\in A_0$. We let 
$$\begin{array}{rlrl}\mu_1&=(p a- p+s,p b +2
p-r-s-3),
 &\mu_2 &=(p a-r-2,p b+r+s+1),\\
\mu_3&=(p a-p+r,p b-p+s), &\mu_4&=\lambda,\\
\mu_5&=(pa -r-s-3,p b+r), 
&\mu_6&=(pa +2 p-r-s-3,p b-p+r),\\ 
\mu_7&=(p a+s,p b-r-s-3), &\mu_8&=(p a+r,p b+s) \\
\mu_9&=(p a+r+s+1,p b-s-2).
\end{array}$$
These weights are depicted in Figure~1~(b). 
Then for $a$ and $b \equiv -1
\pmod p$, $\nabla(\lambda)$ has filtration 
$$\xymatrix@R=10pt@C=5pt
{
& &{\nabla_p(\mu_3)} \ar@{-}[dl] \ar@{-}[dr]  \\
& {\nabla_p(\mu_2)} \ar@{-}[ddr]\ar@{-}[ddrrr]\ar@{-}[dl] 
& &{\nabla_p(\mu_9)} \ar@{-}[dr] \ar@{-}[ddl]\ar@{-}[ddlll] \\
{\nabla_p(\mu_1)} \ar@{-}[d]
& & & &{\nabla_p(\mu_6)} \ar@{-}[d] \\
{{\nabla_p(\mu_7)}} \ar@{-}[drr]   
& & {\nabla_p(\mu_8)} \ar@{-}[d] 
& & {\nabla_p(\mu_5)} \ar@{-}[dll] \\
& &{\nabla_p(\mu_4).}}
$$
For $a \not\equiv -1\pmod p$ there is no extension of
$\nabla_p(\mu_5)$ by $\nabla_p(\mu_6)$.
For $b \not\equiv -1\pmod p$ there is no extension of
$\nabla_p(\mu_7)$ by $\nabla_p(\mu_1)$.
So for $a$ and $b \not\equiv -1\pmod p$ we have:
$$\xymatrix@R=10pt@C=5pt
{
& & &{\nabla_p(\mu_3)} \ar@{-}[dl] \ar@{-}[dr]  \\
& & {\nabla_p(\mu_2)} \ar@{-}[dr]\ar@{-}[drrr]\ar@{-}[dl]\ar@{-}[dll]  
& &{\nabla_p(\mu_9)} \ar@{-}[dr] \ar@{-}[dl]\ar@{-}[dlll]\ar@{-}[drr]  \\
{\nabla_p(\mu_1)} \ar@{-}[drrr]
&{{\ \nabla_p(\mu_7)\!\!\!}} \ar@{-}[drr]   
& &{\nabla_p(\mu_8)} \ar@{-}[d] 
& &{\!\!\!\nabla_p(\mu_5)\ } \ar@{-}[dll] 
&{\nabla_p(\mu_6)} \ar@{-}[dlll] \\
& & &{\nabla_p(\mu_4)}}
$$
and similarly for the other cases for $a$ and $b$.
\end{enumerate}
\end{thm}

\begin{proof}
%
In~\cite{irv}, the structures of certain $\uone$-$B$-modules
$\check{Z}(\lambda)$
are calculated.  We invert these structures to get the
corresponding structure diagrams for the $G_1B$-modules
$\Zone(\lambda)= \Ind_B^{G_1B}(k_\lambda)$ 
as defined in~\cite[II, 9.1]{jantz}.
We will apply~\cite[II, proposition 9.11]{jantz} to these 
diagrams to produce a $p$-filtration of $\nabla(\lambda)$.  
By the remarks of Jantzen~\cite[3.13]{jandar} such filtrations exist
for all $\nabla(\lambda)$ with $\lambda$ dominant and can be obtained
in this way.
 
In all these filtrations there can only be one $\nabla_p(\mu)$
occurring as a sub-module of $\nabla(\lambda)$, namely
$\nabla_p(\lambda)$, as these both have the same simple $G$-socle
$L(\lambda)$. Any module $\nabla_p(\nu)$ appearing above
$\nabla_p(\lambda)$ must have an extension with $\nabla_p(\lambda)$
and this extension has simple socle.
There also can be only one $\nabla_p(\mu)$ which occurs
at the top of $\nabla(\lambda)$, as $\nabla(\lambda)$ has simple
$G$-head by Proposition~\ref{propn:simpsoc}.
This module $\nabla_p(\mu)$ then must have an extension with all
$\nabla_p(\nu)$ appearing below it and this extension must have simple
head.

Case (i). This is clear as $L(p-1,p-1)$ is the Steinberg module.

Case (iii) (and dually case (ii)).
By inverting the diagram in~\cite[theorem 5.2]{irv} and using the 
translation principle 
we have the $G_1B$ composition series for $\Zone(\lambda)$, 
which by applying~\cite[II, proposition 9.11]{jantz} gives us
the $p$-filtration of $\nabla(\lambda)$ as in the statement of the
Lemma above for $b\equiv -1 \pmod p$.
Lemma~\ref{lem:exta} gives us 
the required filtrations when $b \not\equiv -1 \pmod p$.

We now show that the structure obtained so far does not refine further
for $b\equiv -1 \pmod p$.
Consider 
$$\begin{array}{l}
\Ext^1\bigl(\nabla(a-1,b+1)\frob\otimes L(r,s), \nabla(a,b)\frob\otimes
L(s,p-1)\bigr) \\
\cong \Hom_G\bigl(\nabla(a-1,b+1),\nabla(a,b)\otimes\nabla(1,0)\bigr)\\
\cong \Hom_G\bigl(\nabla(a-1,b+1)\otimes\nabla(0,1),\nabla(a,b)\bigr)\\
\cong k
\end{array}$$
where the first isomorphism follows by Lemma~\ref{lem:extone} and
Proposition~\ref{propn:yehia}. 
The last isomorphism follows as
$\nabla(a-1,b+1)\otimes\nabla(0,1)$ has good filtration 
$$\xymatrix@R=15pt{
*=0{\bullet} \ar@{-}[d]^{\textstyle{\nabla(a-1,b+2)}}\\
*=0{\bullet} \ar@{-}[d]^{\textstyle{\nabla(a,b)}}\\
*=0{\bullet} \ar@{-}[d]^{\textstyle{\nabla(a-2,b+1)}} \\
*=0{\bullet} 
}$$
by Lemma~\ref{lem:good} and $(a,b)$ and $(a-1,b+2)$ are not linked 
for $b\equiv -1 \pmod p$.
Let $E$ be the unique (non-split) extension represented by the $\Ext$
group above. We show that $E$ does not
have simple $G$-socle and so does not appear in
$\nabla(pa+s,pb+p-1)$.
From the long exact sequence associated to the short exact sequence for $E$ we have
$$\begin{array}{l}
0\rarr \Hom_{G_1}\bigl(L(r,s),E\bigr)\rarr 
\Hom_{G_1}\bigl(L(r,s), \nabla(a-1,b+1)\frob \otimes L(r,s)\bigr)\\
\rarr\Ext^1_{G_1}\bigl(L(r,s), \nabla(a,b)\frob\otimes L(s,p-1)\bigr) \rarr
\Ext^1_{G_1}\bigl(L(r,s),E\bigr)\rarr 0
\end{array}$$
As 
$$\begin{array}{l}
\Ext^1_{G_1}\bigl(L(r,s), \nabla(a-1,b+1)\frob \otimes L(r,s)\bigr)\\
\cong \nabla(a-1,b+1)\frob \otimes
\Ext^1_{G_1}\bigl(L(r,s),L(r,s)\bigr)\cong 0
\end{array}$$
by Lemma~\ref{lem:extone} and Proposition~\ref{propn:yehia},
and 
$$\Hom_{G_1}\bigl(L(r,s),\nabla(a,b)\frob \otimes L(s,p-1)\bigr)\cong 0$$
We also have
$$\Ext^1_{G_1}\bigl(L(r,s), \nabla(a,b)\frob\otimes L(s,p-1)\bigr)
\cong
\nabla(a,b)\frob\otimes \Ext_{G_1}\bigl(L(r,s),L(s,p-1)\bigr)
\cong \nabla(a,b)\frob\otimes \nabla(1,0)\frob 
$$
by Proposition~\ref{propn:yehia}.

Hence we have an exact sequence
$$\begin{array}{l}
0\rarr \Hom_{G_1}\bigl(L(r,s),E\bigr)\rarr 
\nabla(a-1,b+1)\frob \\
\hspace{40pt}
\stackrel\phi\rarr \nabla(a,b)\frob\otimes
\nabla(1,0)\frob
\rarr \Ext^1_{G_1}\bigl(L(r,s),E\bigr)\rarr 0
\end{array}$$
If $\phi$ is injective then $\Hom_{G_1}\bigl(L(r,s),E\bigr)\cong 0$. 
We claim that this is not the case.

Now 
$$\begin{array}{rl}
\phi\mbox{ injective}&\impl \nabla(a-1,b+1)\frob\mbox{ embeds in } 
\nabla(a,b)\frob \otimes \nabla(1,0)\frob\\
 &\lrimpl  \nabla(a-1,b+1)\mbox{ embeds in }
\nabla(a,b) \otimes \nabla(1,0)
\end{array}$$
but $\nabla(a-1,b+1)$ cannot embed in
$\nabla(a,b)\otimes\nabla(1,0)$
as
$$\begin{array}{rl}
k &\cong \Hom_G\bigl(\nabla(a-1,b+1),\nabla(a,b-1)\bigr)\\
& \subseteq
\Hom_G\bigl(\nabla(a-1,b+1),\nabla(a,b)\otimes\nabla(1,0)\bigr)\\
& \cong
\Hom_G\bigl(\nabla(a-1,b+1)\otimes\nabla(0,1),\nabla(a,b)\bigr)
\cong k,
\end{array}$$
where the first isomorphism follows as $(a-1,b+1)$ and $(a,b-1)$
differ by a single reflection~\cite[II, corollary 6.24]{jantz}  and the last $\Hom$ group we calculated above.
However homomorphisms in the first $\Hom$
group are clearly not 1-1. Hence $\phi$ cannot be injective.
Thus we have $\Hom_{G_1}\bigl(L(r,s),E\bigr)\not\cong 0$. 

Now using~\cite[2.2 (1)]{jandar}
$$\begin{array}{rl}
\soc_G(E) &=\soc_G\Bigl(\Hom_{G_1}\bigl(L(s,p-1),E\bigr)\Bigr)\otimes
L(s,p-1)\\
&\hspace{40pt}
\oplus \soc_G\Bigl(\Hom_{G_1}\bigl(L(r,s),E\bigr)\Bigr)\otimes L(r,s)\\
&=\soc_G\bigl(\nabla(a,b)\frob \otimes L(s,p-1)\bigr)\\
&\hspace{40pt}
\oplus \soc_G\Bigl(\Hom_{G_1}\bigl(L(r,s),E\bigr)\Bigr)\otimes L(r,s)
\end{array}$$
Both terms are non-zero. Hence $E$ does not have simple $G$-socle.

Case (iv). We write $\lambda=p(a,b)+(r,s)$ with $r+s=p-2$, and $(a,b)\in X^+$.
By inverting the diagram in~\cite[theorem 3.3]{irv} and using 
the translation principle 
we see that $\Zone(\lambda)$
has $G_1T$ composition series:
$$\xymatrix@R=10pt@C=-40pt
{
& {\Lone\bigl(p(a-1,b-1)+(r,s)\bigr)}\ar@{-}[dl] \ar@{-}[dr]
\\
{\Lone\bigl(p(a,b-1)+(p-1,r)\bigr)} \ar@{-}[dr] &
 &{\Lone\bigl(p(a-1,b)+(s,p-1)\bigr)} \ar@{-}[dl] \\
&{\Lone\bigl(p(a,b)+(r,s)\bigr).}}  
$$
So we need to show that
$$\begin{array}{l}
\Ext^1\bigl(\nabla(a,b-1)\frob\otimes
L(p-1,r),\nabla(a-1,b)\frob\otimes L(s,p-1)\bigr)\\
\cong \Ext^1\bigl(\nabla(a-1,b)\frob\otimes L(s,p-1),
\nabla(a,b-1)\frob \otimes L(p-1,r)\bigr)\\
\cong 0.
\end{array}$$
But $\Ext_{G_1}^1\bigl(L(s,p-1),L(p-1,r)\bigr)\cong
\Ext_{G_1}^1\bigl(L(p-1,r),L(s,p-1)\bigr)\cong 0$
by Proposition~\ref{propn:yehia}, and  so
by Lemma~\ref{lem:extone} we are done. 

Case (v) (and dually case (vi)). 
We know that this filtration is correct on the level of characters
using~\cite[theorem 3.1]{jandar} and the results of~\cite{irv}. 
We also know that $\nabla_p(\lambda)$ embeds in
$\nabla(\lambda)$. Further, for $p\ne 3$ 
$\Ext^1_{G_1}\bigl(L(s,p-r-s-3),L(r,s)\bigr)=0$ by Proposition~\ref{propn:yehia}
so by Lemma~\ref{lem:extone} we have
$$\Ext_G^1\bigl(\nabla(a-2,0)\frob\otimes
L(s,p-r-s-3),\nabla(a,0)\frob\otimes L(r,s)\bigr)\cong 0.$$
For $p=3$, $(s,p-r-s-3)=(r,s)=(0,0)$ and we then have
$$\begin{array}{l}
\Ext_G^1\bigl(\nabla(a-2,0)\frob\otimes
L(s,p-r-s-3),\nabla(a,0)\frob\otimes L(r,s)\bigr) 
\\ 
\cong \Ext_G^1\bigl(\nabla(a-2,0),\nabla(a,0)\bigr) 
\cong 0
\end{array}$$
using Lemma~\ref{lem:extone} and observing that  $(a-2,0)$ is not 
comparable to $(a,0)$ in the usual ordering of $G$-weights.
Hence we get the
filtration as in the statement of the Proposition.

Case (viii). 
By inverting the diagram in~\cite[theorem 5.3]{irv}
we get the $G_1B$ structure of $\Zone(\lambda)$ with
$\lambda = (p-s-2,p-r-2)+p(a,b)$ and
$(r,s)\in A_0$.
So applying~\cite[II, proposition 9.11]{jantz}
we get a filtration as stated for the
case $a$ and $b \equiv -1 \pmod p$.

If $p\ne 3$, then  by
Lemma~\ref{lem:extzero} none of the $\nabla_p(\mu_i)$ which appear 
above $\nabla_p(\lambda)$ can have non-split extensions by $\nabla_p(\mu_3)$.
Also $\nabla_p(\lambda)$ does not have an extension by
$\nabla_p(\mu_2)$ nor by $\nabla_p(\mu_9)$.
Now consider $p=3$. We write $\mu_i=p\mu_{i1} +\mu_{i0}$.
For $i=1$, $5$ , $6$, $7$ and  $8$ Lemma~\ref{lem:extone} gives us
$$\Ext^1\bigl(\nabla_p(\mu_3),\nabla_p(\mu_i)\bigr)
\cong \Ext^1\bigl(\nabla(\mu_{31}), \nabla(\mu_{i1})\bigr)
\cong 0$$
as none of the $\mu_{i1}$ are less than $\mu_{31}$.

We now take $p \ge 3$ again. We need to show that the filtration
simplifies for $a$ (or $b$) 
$\not \equiv -1 \pmod p$, and that it doesn't simplify for $a$ (or
$b$) $\equiv -1 \pmod p$.

Lemmas~\ref{lem:exta} and~\ref{lem:extb} give us
the filtration as stated in the Theorem, part (viii) when $a$ (or $b$) 
$\not \equiv -1\pmod p$. 

Now the same argument as for case (iii) (with $p \ge 5$)
shows that both the non-split extensions
$E_1$ and $E_2$ defined via:
$$\begin{array}{c}
0 \rarr\nabla_p(\lambda)\rarr E_1 \rarr \nabla_p(\mu_1)\rarr 0\\
0 \rarr\nabla_p(\lambda)\rarr E_2 \rarr \nabla_p(\mu_6)\rarr 0
\end{array}$$
do not have simple socle and so the filtration does not refine any
further for $a$ (or $b$) $\equiv  -1\pmod p$.
If $p=3$ then the argument used in case (iii) still carries through,
as when we remove the Frobenius twist then block considerations also
alllow us to remove the other direct summands.

Case (vii). We get the $\uone$-$B$-filtration of $\check{Z}(\lambda)$ by
taking its dual and then applying~\cite[theorem 5.3]{irv}. We then
invert this structure to get the $G_1B$ structure of $\Zone(\lambda)$.
Then applying~\cite[II, proposition 9.11]{jantz} gives us
the corresponding filtration of $\nabla(\lambda)$.     
The same argument as before shows that there is no non-split extension of
$\nabla_p(\lambda)$ by $\nabla_p(\mu_i)$ for $i=3$, $5$, $6$, $7$ and
$8$, and also that there is no non-split extension of $\nabla_p(\mu_i)$ by
$\nabla_p(\mu_9)$ for $i=2$ and $4$.
The extensions of $\nabla_p(\mu_5)$ by $\nabla_p(\mu_6)$ and of
$\nabla_p(\mu_3)$ by $\nabla_p(\mu_8)$ are split as in Case (vii) for 
$a$ or $b\not \equiv 0\pmod p$ using Lemmas~\ref{lem:exta} and~\ref{lem:extb}.

Now suppose $a\equiv 0\pmod p$. (The case with $b\equiv 0\pmod p$ is
similar.)
Consider the extension $E$ defined via the short exact sequence
$$0\rarr \nabla_p(\mu_5)\rarr E\rarr \nabla_p(\mu_9)\rarr 0.$$
Such a non-split extension $E$, exists and is unique, since,
$$\Ext^1_G\bigl(\nabla_p(\mu_9),\nabla_p(\mu_5)\bigr)
\cong \Hom_G\bigl(\nabla(a-1,b-1),\nabla(a-2,b)\otimes
\nabla(0,1)\bigr)$$
using Lemma~\ref{lem:extone} and Proposition~\ref{propn:yehia}.
But this is isomorphic to $k$, as the costandard modules
($\nabla(a-2,b+1)$, $\nabla(a-1,b-1)$ and $\nabla(a-3,b)$) appearing in a good
filtration of $\nabla(a-2,b)\otimes
\nabla(0,1)$ are not linked for $a\equiv 0\pmod p$.
We have a long exact sequence:
$$\begin{array}{rl}
0&\rarr \Hom_{G_1}\bigl(\nabla(\mu_9),L(s,p-r-s-3)\bigr)
\rarr \Hom_{G_1}\bigl(E,L(s,p-r-s-3)\bigr)\hspace{10pt}
\\
&\rarr \Hom_{G_1}\bigl(\nabla(\mu_5),L(s,p-r-s-3)\bigr)
\rarr \Ext^1_{G_1}\bigl(\nabla(\mu_9),L(s,p-r-s-3)\bigr)
\\
&\rarr \Ext^1_{G_1}\bigl(E,L(s,p-r-s-3)\bigr)
\rarr \Ext^1_{G_1}\bigl(\nabla(\mu_5),L(s,p-r-s-3)\bigr).
\end{array}$$
The first $\Hom$ group and the last $\Ext^1$ group are both $0$.
The third $\Hom$ group is known and the first $\Ext^1$ group is known 
using Proposition~\ref{propn:yehia}.
So we have, for $p\ge 5$,
(for $p=3$ the $\nabla(0,1)\frob$ would be replaced by $k \oplus
\nabla(1,0)\frob \oplus \nabla(0,1)\frob$):
$$\begin{array}{l}
0\rarr \Hom_{G_1}\bigl(E,L(s,p-r-s-3)\bigr)
\rarr\Delta(b,a-2)\frob\\
\hspace{40pt} \stackrel\phi\rarr 
\Delta(b-1,a-1)\frob\otimes \nabla(0,1)\frob
\rarr \Ext^1_{G_1}\bigl(E,L(s,p-r-s-3)\bigr)
\end{array}$$
We show that $\phi$ is not injective.
Now 
$$\begin{array}{rl}\phi \in 
&\Hom_{G/{G_1}}\bigl(\Delta(b,a-2)\frob,
\Delta(b-1,a-1)\frob\otimes\nabla(0,1)\frob\bigr)\\
&\cong \Hom_G\bigl(\Delta(b,a-2),
\Delta(b-1,a-1)\otimes\Delta(0,1)\bigr)\\
&\cong \Hom_G\bigl(\nabla(a-1,b-1),
\nabla(a-2,b)\otimes\nabla(0,1)\bigr)\cong k.
\end{array}$$

Also, 
$$\begin{array}{rl}
\phi\mbox{ injective}&\lrimpl \Delta(b,a-2)\frob\mbox{ embeds in } 
\Delta(b-1,a-1)\frob \otimes \Delta(0,1)\frob\\
&\lrimpl  \Delta(b,a-2)\mbox{ embeds in }
\Delta(b-1,a-1) \otimes \Delta(0,1)\\
&\lrimpl  \nabla(a-2,b)\mbox{ is a quotient of }
\nabla(a-1,b-1) \otimes \nabla(1,0).
\end{array}$$
However, $\nabla(a-2,b)$ cannot be quotient of
$\nabla(a-1,b-1)\otimes\nabla(1,0)$
as
$$\begin{array}{rl}k 
& \cong \Hom_G\bigl(\nabla(a,b-1),\nabla(a-2,b)\bigr)\\
&\subseteq
\Hom_G\bigl(\nabla(a-1,b-1)\otimes
\nabla(1,0),\nabla(a-2,b)\bigr)\\
&\cong \Hom_G\bigl(\nabla(a-1,b-1),
\nabla(a-2,b)\otimes\nabla(0,1)\bigr)
\cong k,
\end{array}$$
where the first isomorphism follows as $(a,b-1)$ and $(a-2,b)$
differ by a single reflection~\cite[II, corollary 6.24]{jantz}. 
But homomorphisms in the first $\Hom$
group are clearly not onto. Hence $\phi$ cannot be injective.
and so we have that $\Hom_{G_1}\bigl(E, L(s,p-r-s-3)\bigr)$ is non-zero.
By a similar argument to that before, (using the head functor $\hd$,
in place of the socle functor $\soc$),
we have that $E$ does not have simple $G$-head and so it does not
appear in $\nabla(\lambda)$ which has simple
head using Lemma~\ref{lem:weylsoc}.
The argument for $p=3$ is essentially the same. As soon as we remove
the Frobenius twists ($\frob$) from the modules then block
considerations allow us to remove the extra summands that appear.
\end{proof}

\section{Good filtration dimensions and global dimensions for $\SL_3$.}

In this section we calculate the good filtration dimension 
of $\nabla_p(\lambda)$. There are two cases
to consider, one where $\lambda$ is on a wall and the other where $\lambda$ is
inside an alcove. We consider this latter case first.

\begin{lem}\label{lem:mzero}
Let $\lambda = p(a,b)+(r,s)$ with $(r,s)\in A_0$, $(a,b)\in X^+$
and $b\ge 1$. 
We take $M_{\lambda}$ be the submodule of $\nabla(\lambda)$ with 
$p$-filtration labelled by $\mu_1$ through to $\mu_3$ in the diagram
of Proposition~\ref{thm:pfiltra} part (vii). (If $b=1$ we take $\nabla(a,b-2)=0$.) 
Then $M_\lambda$ has good resolution
$$0\rarr M_\lambda \rarr \nabla(\lambda_0) \rarr \nabla(\lambda_1)\rarr \cdots \rarr
\nabla(\lambda_a)\rarr 0$$
where $\lambda_{2i}=p(a-2i,b+i)+(r,s)$ and
$\lambda_{2i+1}=p(a-2i-1,b+i)+(p-r-2,r+s+1)$
for $0\le i\le \lfloor \frac{a}{2}\rfloor$.
%

Similarly, let $\mu= p(a,b) +(p-r-2,r+s+1)$ with $(p-r-2,r+s+1)$ in the upper
alcove of $X_1$ and $(a,b)\in X^+$.
We take $M_{\mu}$ be the submodule of
$\nabla(\mu)$ with
$p$-filtration labelled by $\mu_4$ through to $\mu_9$ in the diagram
of Proposition~\ref{thm:pfiltra} part (viii).
Then $M_\mu$ has good resolution 
$$0\rarr M_\mu \rarr \nabla(\mu_0) \rarr \nabla(\mu_1)\rarr \cdots \rarr
\nabla(\mu_a)\rarr 0$$
where $\mu_{2i}=p(a-2i,b+i)+(p-r-2,r+s+1)$ and
$\mu_{2i+1}=p(a-2i-1,b+i+1)+(r,s)$
for $0\le i\le \lfloor \frac{a}{2}\rfloor$.
\end{lem}
\begin{proof} We proceed by induction on a.
If $a=0$ then $M_\lambda = \nabla(\lambda)$ and $M_\mu= \nabla(\mu)$
by Proposition~\ref{thm:pfiltra} and we are done.
For $a \ge 1$ it is sufficient to show the existence of two short exact
sequences
$$\begin{array}{l} 
0 \rarr M_{\lambda} \rarr \nabla(\lambda) \rarr M_{\lambda_1} \rarr 0 \\
0 \rarr M_{\mu} \rarr \nabla(\mu) \rarr M_{\mu_1} \rarr 0. 
\end{array}$$
Now $\Hom\bigl(\nabla(\lambda), \nabla(\lambda_1)\bigr) \cong k$, as
$\lambda$ and $\lambda_1$ satisfy the conditions in~\cite[II, corollary
6.24]{jantz}. So we need to show that $M_{\lambda_1}$ is the image of
this unique homomorphism $\phi$ (as suggested by the labelling of the $\mu$'s
in Proposition~\ref{thm:pfiltra}) and that $M_\lambda$ is contained
in the kernel of $\phi$.
Now $\nabla(\lambda_1)$ has simple socle $L(\lambda_1)$ 
which appears just once right above $M_\lambda$ in $\nabla(\lambda)$.
(As all $\nabla_p(\mu)$ which occur above $\nabla_p(\lambda_1)$
have $\mu < \lambda_1$.)
Hence $M_{\lambda_1}$ is the image of $\phi$ and we are done.
The second sequence is similar.
\end{proof}

\begin{rem}
We have identified the images of the maps in the resolution (1)
of Jantzen~\cite[3.12, remark 2]{jandar}.
\end{rem}

\begin{cor}\label{cor:gfdmzero}
For $M_{\lambda}$ and $M_{\mu}$ as defined above we have
$\gfd(M_{\lambda})=\gfd(M_{\mu})=a$.
Furthermore for $\tau \in X^+$ we have
$$\Ext^a\bigl(\Delta(\tau), M_\lambda \bigr) \cong 
\Hom \bigl(\Delta(\tau), \nabla(\lambda_a)\bigr)\cong \delta_{\tau
\lambda_a}k$$
and
$$\Ext^a\bigl(\Delta(\tau), M_\mu \bigr) \cong 
\Hom \bigl(\Delta(\tau), \nabla(\mu_a)\bigr)\cong \delta_{\tau \mu_a}k$$
\end{cor}
\begin{proof}
The result follows by dimension shifting and by noting that
$\Hom\bigl(\Delta(\tau),\nabla(\lambda_{a-1})\bigr)$ is zero if
$\Hom\bigl(\Delta(\tau),\nabla(\lambda_{a})\bigr)$ is non-zero.
The argument for $M_\mu$ is similar.
\end{proof}

Suppose $\lambda = p(a,b) + \nu$ and $\nu \ne (p-1, p-1)$.
We define $g(\lambda)$ as follows
$$g(\lambda)=  \left\{ \begin{array}{ll}
                                2(a+b)   &\mbox{if $\nu$ is inside 
                                         a lower alcove}\cr
                                2(a+b)+1 &\mbox{if $\nu$ is inside 
                                         an upper alcove}\cr
                                a+b      & \mbox{if $\nu$ lies on a
                                         wall.}
                                \end{array}
\right. $$
\noindent We will eventually show that $\gfd\bigl(L(\lambda)\bigr)=g(\lambda)$.

\begin{lem}\label{lem:gdecr}
If $\mu < \lambda$ and $\mu \in \B(\lambda)$ then $g(\mu)\le g(\lambda)$.
\end{lem}
\begin{proof}
Suppose $\mu =(c,d) < \lambda =(a,b)$.
If we think of $\lambda$ as the $S(3,a+2b)$ weight $(a+b,b,0)$ then 
$\mu$ can be thought of the $S(3,a+2b)$ weight $(c+d+e,d+e,e)$ for an
appropriate $e \in \N$.
Since $\mu < \lambda$ we have $(c+d+e,d+e,e) < (a+b,b,0)$ in the 
dominance ordering on partitions.
But this corresponds to `falling boxes' in the corresponding Young diagram.
Consequently it is clear that $g(\lambda)$ cannot be increased by
moving boxes further down the diagram, hence the result.
\end{proof}

\begin{lem}\label{lem:Qweight3}
Suppose $\lambda$ is primitive and define $Q$ to be the quotient
$\nabla_p(\lambda)/L(\lambda)$.
Then $$g\bigl(\hw(Q)\bigr) \le g(\lambda)-1.$$
\end{lem}
\begin{proof}
We write $\lambda=p(a,b)+\nu$. We define $\bar{g}(\lambda)=a+b$.
By assumption $\nu \ne (p-1,p-1)$.
We have  $L(\lambda) \cong L(a,b)\frob \otimes L(\nu)$
by Steinberg's tensor product theorem and so $Q\cong R\frob\otimes
L(\nu)$ where $R$ is the quotient $\nabla(a,b)/L(a,b)$.
If we write $(a,b)=p(a^\prime,b^\prime) +\nu^\prime$ with $\nu^\prime$
$p$-restricted, 
then it is clear using~\ref{thm:pfiltra} and induction that 
$$\hw(R)= \left\{ \begin{array}{ll}
                        p\hw\bigl(\nabla(a^\prime,b^\prime)/L(a^\prime,
                        b^\prime)\bigr) + \nu^\prime
                                &\mbox{if $\nu^\prime= (p-1,p-1)$}\cr
                        p(a^\prime+1,b^\prime-1)+(r,p-r-2)
                                &\mbox{if $\nu^\prime= (p-1,r)$}\cr
                        p(a^\prime-1,b^\prime+1)+(p-s-2,s)
                                &\mbox{if $\nu^\prime= (s,p-1)$}\cr
                        \{p(a^\prime,b^\prime-1)+(p-1,r),
                              &\mbox{if $\nu^\prime= (r,s)$ with $r+s=p-2$}\cr
                        \ \ \ p(a^\prime-1,b^\prime)+(s,p-1)\}
                                          & \cr
                        \{p(a^\prime,b^\prime-1)+(r+s+1,p-s-2),
                                 &\mbox{if $\nu^\prime= (r,s)\in A_0$}\cr
                        \ \ \ p(a^\prime-1,b^\prime)+(p-r-2,r+s+1)\}
                                 & \cr
                        \{p(a^\prime,b^\prime)+(r,s),
                                 &\mbox{if $\nu^\prime=
                                 (p-s-2,p-r-2)$}\cr
                        \ \ \ p(a^\prime-1,b^\prime+1)+(s,p-r-s-3),
                                 &\mbox{\ \ \ inside an upper alcove.}\cr
                        \ \ \ p(a^\prime+1,b^\prime-1)+(p-r-s-3,r)\}
                                 & 
                \end{array} \right. 
$$
%
Since $\hw(Q)=p(\hw(R))+\nu$ we have 
$$\bar{g}(\hw(Q))= \left\{ \begin{array}{ll}
                        \bar{g}\bigl(p^2\hw\bigl(\nabla(a^\prime,b^\prime)
                             /L(a^\prime, b^\prime)\bigr)\bigr)
                                           +2p-2 
                                &\mbox{if $\nu^\prime= (p-1,p-1)$}\cr
                        p(a^\prime+b^\prime)+p-2
                            &\mbox{if $\nu^\prime= (p-1,r)$
                                 \ or $\nu^\prime= (s,p-1)$}\cr
                        \{p(a^\prime+b^\prime)+r-1,  p(a^\prime+b^\prime)+s-1\}
                              &\mbox{if $\nu^\prime= (r,s)$ with $r+s=p-2$}\cr
                        \{p(a^\prime+b^\prime)+r-1,
                         p(a^\prime+b^\prime)+s-1\}
                                 &\mbox{if $\nu^\prime= (r,s)\in A_0$}\cr
                        \{p(a^\prime+b^\prime)+r+s,
                                 &\mbox{if $\nu^\prime=
                                 (p-s-2,p-r-2)$}\cr
                        \qquad \ \  \ \ p(a^\prime+b^\prime)+p-r-3,
                                 &\mbox{\ \ \ inside an upper alcove.}\cr
                        \qquad \qquad \  \ \ p(a^\prime+b^\prime)+p-s-3\}
                                 & 
                \end{array} \right. 
$$
Now if $\nu^\prime = (c,d)$ then $\bar{g}(\lambda)=
p(a^\prime+b^\prime)+c+d$ and so we have 
$$\bar{g}(\hw(Q))= \left\{ \begin{array}{ll}
                        \bar{g}\bigl(p^2\hw(R^\prime)\bigr) +2p-2 
                                &\mbox{if $\nu^\prime= (p-1,p-1)$}\cr
                        \bar{g}(\lambda) - (r+1) 
                                &\mbox{if $\nu^\prime= (p-1,r)$}\cr
                        \bar{g}(\lambda) - (s+1)
                                &\mbox{if $\nu^\prime= (s,p-1)$}\cr
                        \{\bar{g}(\lambda) - (s+1),\ 
                         \bar{g}(\lambda) - (r+1)\}
                              &\mbox{if $\nu^\prime= (r,s)$ with $r+s=p-2$}\cr
                        \{\bar{g}(\lambda) - (s+1),
                        \ \bar{g}(\lambda) - (r+1)\}
                                 &\mbox{if $\nu^\prime= (r,s)\in A_0$}\cr
                        \{\bar{g}(\lambda) - 2(p-r-s-2),
                                 &\mbox{if $\nu^\prime=
                                 (p-s-2,p-r-2)$}\cr
                        \qquad \ \ \bar{g}(\lambda) - (p-s-1),
                                 &\mbox{\ \ \ inside an upper alcove}\cr
                        \qquad \qquad \ \bar{g}(\lambda) - (p-r-1)\}
                                 & 
                \end{array} \right. 
$$
where $R^\prime = \nabla(a^\prime, b^\prime)/L(a^\prime, b^\prime)$.
We write $(a^\prime, b^\prime)=p(a^{\prime\prime},
b^{\prime\prime})+\nu^{\prime\prime}$. 
Now $\bar{g}\bigl(\hw(R^\prime)\bigr) \le a^{\prime\prime} + b^{\prime\prime}$. 
So if $\nu^\prime =(p-1,p-1)$ we have
$\bar{g}(\hw(Q))\le p^2(a^{\prime\prime}+b^{\prime\prime})+2p-2
<p(a^\prime+b^\prime)+2p-2=g(\lambda)$ provided $\nu^{\prime\prime}\ne
0$.
If $\nu^{\prime\prime}=0$ then $\bar{g}(R^\prime) = a^{\prime\prime} +
b^{\prime\prime}-1$ and so we have
 $\bar{g}(\hw(Q))= p^2(a^{\prime\prime}+b^{\prime\prime})-p^2+2p-2
<p(a^\prime+b^\prime)+2p-2=g(\lambda)$. Thus the result follows for
$\nu^\prime=(p-1,p-1)$.
Since $0 \le r,s \le p-2$ and $r+s < p-2$ inside an alcove the result
follows for all other $\nu^\prime$.
%
\end{proof}

The next Proposition shows that $\gfd\bigl(L(\lambda)\bigr)
=g(\lambda)$ for $\lambda$ inside an alcove.
The following Lemma forms part of the inductive step.

\begin{lem}\label{lem:ind2}
Let $\lambda \in X^+$ and 
suppose $g(\mu)=\gfd\bigl(L(\mu)\bigr)= \gfd(\bigl(\nabla_p(\mu)\bigr)$ for
all $\mu < \lambda$ with $\mu \in \B(\lambda)$ and 
$\gfd\bigl(\nabla_p(\lambda)\bigr)=g(\lambda)$. Then
$\gfd\bigl(L(\lambda)\bigr)=\gfd\bigl(\nabla_p(\lambda)\bigr)$.
\end{lem}
\begin{proof}
We have a short exact sequence
$$0\rarr L(\lambda) \rarr \nabla_p(\lambda) \rarr Q \rarr 0.$$
Lemma~\ref{lem:Qweight3} gives us $\gfd\bigl(\nabla_p(\lambda)\bigr)> \gfd(Q)
$, since all the composition factors of $Q$ have weights with smaller
good filtration dimension,
and the result follows by Lemma~\ref{lem:gandi} part (i).
\end{proof}

\begin{thm}\label{thm:gfdalc}
Let $\nu=(r,s)$ be a weight inside the fundamental alcove and
$\bar{\nu}=(p-s-2,p-r-2)$ its reflection in the upper alcove. Then we have,
$$\begin{array}{rl}
 &\gfd\bigl(L(p(a,b)+\nu)\bigr)
=\gfd\bigl(\nabla(a,b)\frob\otimes L(\nu)\bigr)=2(a+b)
\\
\mbox{and} &\gfd\bigl(L(p(a,b)+\bar{\nu})\bigr)
=\gfd\bigl(\nabla(a,b)\frob\otimes L(\bar{\nu})\bigr)=2(a+b)+1.
\end{array}
$$
Moreover 
$$\begin{array}{rl}
 &\Ext^{4(a+b)}\bigl(\Delta(a,b)\frob\otimes
L(\nu),\nabla(a,b)\frob\otimes L(\nu)\bigr)\not
\cong 0 \\
\mbox{and} &
\Ext^{4(a+b)+2}\bigl(\Delta(a,b)\frob\otimes L(\bar{\nu}),
\nabla(a,b)\frob\otimes L(\bar{\nu})\bigr) \not\cong  0.
\end{array}$$
\end{thm}

\begin{proof} We proceed by induction on $a+b$.
For $a+b=0$ we have $L(\nu)=\nabla(\nu)=\Delta(\nu)$ and 
$$\gfd\bigl(L(\nu)\bigr)=\gfd\bigl(\nabla(\nu)\bigr)=0.$$ 
Also we have a non-split short exact sequence for $L(\bar{\nu})$, namely
$$0 \rarr L(\bar{\nu}) \rarr \nabla(\bar{\nu})\rarr L(\nu) \rarr 0.$$
Hence we have $\gfd\bigl(L(\bar{\nu})\bigr)=1$ and
$\Ext^2\bigl(L(\bar{\nu}),L(\bar{\nu})\bigr) \not\cong 0$ by 
Lemma~\ref{lem:globtwo}. 

Now suppose $a+b\ge 1$.
Let $\lambda = p(a,b)+\nu$ and $M_\lambda$ be as in Lemma~\ref{lem:mzero}.
We first show $\gfd\bigl(\nabla_p(\lambda)\bigr)=g(\lambda)$.

Case (i): Suppose $b\ge 2$. We have a short exact sequence
\begin{equation}\label{Mses1}
0\rarr \nabla(a,b)\frob\otimes L(r,s)\rarr M_\lambda \rarr Q\rarr 0
\end{equation}
where $Q$ is the (unique) extension
\begin{equation}\label{Qses1}
0\rarr \nabla(a,b-1)\frob\otimes L(r+s+1,p-s-2)
 \rarr Q
\rarr\nabla(a,b-2)\frob\otimes L(p-r-s-3,r)\rarr 0.
\end{equation}
By 
induction we have $\gfd\bigl(\nabla(a,b-1)\frob\otimes
L(r+s+1,p-s-2)\bigr)=2(a+b)-1$ and
$\gfd\bigl(\nabla(a,b-2)\frob\otimes L(p-r-s-3,r)\bigr)=2(a+b)-4$.
Thus we have $\gfd(Q)=2(a+b)-1$ by Lemma~\ref{lem:gandi}. We also have
$\gfd(Q)>\gfd(M_\lambda)=a$ by Corollary~\ref{cor:gfdmzero}, so we have
$\gfd\bigl(\nabla(\lambda)\bigr)=2(a+b)$ by Lemma~\ref{lem:gandi}, as required.
Lemma~\ref{lem:gandi} also gives us
$$\begin{array}{l}\Ext^{4(a+b)}\bigl(\Delta_p(\lambda),\nabla_p(\lambda)\bigr)
\cong \Ext^{4(a+b)-1}\bigl(\Delta_p(\lambda),Q\bigr)
\\ \hspace{10pt}
\cong \Ext^{4(a+b)-1}\bigl(\Delta_p(\lambda),\nabla(a,b-1)\frob\otimes
L(r+s+1,p-s-2)\bigr)
\\ \hspace{10pt}
\cong \Ext^{4(a+b)-2}\bigl(Q^\circ, \nabla(a,b-1)\frob\otimes
L(r+s+1,p-s-2)\bigr)
\\ \hspace{10pt}
\cong \Ext^{4(a+b)-2}\bigl(\Delta(a,b-1)\frob\otimes
L(r+s+1,p-s-2), 
\nabla(a,b-1)\frob\otimes
L(r+s+1,p-s-2)\bigr)
\end{array}
$$
where the first isomorphism follows using sequence~(\ref{Mses1}), the
second from~(\ref{Qses1}), the third from the dual of~(\ref{Mses1}) and the
fourth from the dual of~(\ref{Qses1}).
This last $\Ext$ group is non-zero by induction.

Case (ii): Suppose $b=1$. The argument above simplifies as $Q=
\nabla(a,0)\frob\otimes L(r+s+1,p-s-2)$ and so by induction we
have $\gfd(Q)= 2(a+1)-1>\gfd(M_\lambda)=a$. Hence
$\gfd\bigl(\nabla_p(\lambda)\bigr)=2(a+1)$
and
$$\Ext^{4(a+1)}\bigl(\Delta_p(\lambda),\nabla_p(\lambda)\bigr)\not\cong 0$$
as required.

Case (iii): Suppose $b=0$. We have a short exact sequence
$$0\rarr \nabla(a,0)\frob\otimes L(r,s)\rarr \nabla(\lambda) \rarr Q\rarr 0$$ 
where $Q$ is the (unique) extension
$$
0\rarr \nabla(a-1,0)\frob\otimes L(p-r-2,r+s+1)
\rarr Q
\rarr
\nabla(a-2,0)\frob\otimes L(s,p-r-s-3)\rarr 0.
$$
By induction we have $\gfd\bigl(\nabla(a-1,0)\frob\otimes
L(r+s+1,p-s-2)\bigr)=2a-1$ and
$\gfd\bigl(\nabla(a-2,0)\frob\otimes L(s, p-r-s-3)\bigr)=2a-4$.
A similar argument to Case (i) yields the required result.


\vskip5pt
Now let $\mu = p(a,b)+\bar{\nu}$ and $M_\mu$ be as in Lemma~\ref{lem:mzero}.
We have a short exact sequence
$$0 \rarr \nabla_p(\mu) \rarr M_\mu \rarr Q \rarr 0.$$
Thus $\gfd\bigl(\nabla_p(\mu)\bigr)=\gfd(Q)+1$ provided $\gfd(Q)> \gfd(M_\mu)
=a$.
We have that $\gfd(Q) \le 2(a+b)$ by induction, so $\gfd(\nabla_p(\mu))\le 2(a+b)+1$.
We need to show that both these bounds are attained.

Suppose $b\ge 1$. Define $R$ via the short exact sequence
$$0\rarr R \rarr Q \rarr \nabla(a,b-1)\frob\otimes L(r+s+1,p-s-2)
\rarr 0.$$

Here (using the notation of Proposition~\ref{thm:pfiltra})
it is clear by induction and using the $p$-filtration of $R$ 
that for $M$ a $G$-module we have
$$
\Ext^{\wfd(M)+2(a+b)}(M,R)
\cong \Ext^{\wfd(M) + 2(a+b)}\bigl(M,\nabla_p(\mu_6)\bigr)
\oplus \Ext^{\wfd(M)+2(a+b)}\bigl(M,\nabla_p(\mu_8)\bigr)
$$
since the other $\nabla_p(\mu_i)$ that appear in $R$ have good
filtration dimension equal to $2(a+b)-2$ and so cannot contribute to
this $\Ext$ group. We have a direct sum since there is no extension
appearing between $\nabla_p(\mu_6)$ and $\nabla_p(\mu_8)$.

The long exact sequence gives us
$$\begin{array}{l}
\Ext^{\wfd(M) + 2(a+b)-1}\bigl(M,\nabla(a,b-1)\frob \otimes
L(r+s+1,p-s-2)\bigr) \\
\hspace{55pt}\rarr \Ext^{\wfd(M)+2(a+b)}\bigl(M, R\bigr)
\rarr  \Ext^{\wfd(M)+2(a+b)}(M, Q) \rarr 0.
\end{array}$$
But the middle $\Ext$ group is as above.
%
Also 
$$\begin{array}{l}
\Ext^{\wfd(M)+2(a+b)-1}\bigl(M,\nabla(a,b-1)\frob \otimes
L(r+s+1,p-s-2)\bigr) \\
\hspace{100pt}\cong \Ext^{\wfd(M)+2(a+b)}\bigl(M,\nabla(a,b)\frob \otimes
L(r,s))\bigr)
\end{array}$$
using the case above.
Thus we have
$$\begin{array}{l}\Ext^{\wfd(M)+2(a+b)}\bigl(M,
        \nabla(a,b)\frob \otimes L(r,s)\bigr) \\
\hspace{80pt}
\rarr \Ext^{\wfd(M)+2(a+b)}\bigl(M,
        \nabla(a,b)\frob \otimes L(r,s)\bigr) \\
 \hspace{100pt}
\oplus \Ext^{\wfd(M)+2(a+b)}\bigl(M,
        \nabla(a+1,b-1)\frob\otimes L(p-r-s-3,r)\bigr) \\
\hspace{220pt}
\rarr \Ext^{\wfd(M)+2(a+b)}(M, Q) \rarr 0.
\end{array}$$
Hence $\Ext^{\wfd(M)+2(a+b)}(M, Q)$ will be non-zero 
if $$\Ext^{\wfd(M)+2(a+b)}\bigl(M,
\nabla(a+1,b-1)\frob\otimes L(p-r-s-3, r)\bigr) \ne 0.$$

We know by induction that
$\gfd \bigl(\nabla(a+1,b-1)\frob\otimes L(p-r-s-3,r)\bigr)
= 2(a+b)$ and $\gfd(Q)\le 2(a+b)$.
If we take $M\in\doog$  (so $\wfd(M)=0$)
with the last $\Ext$ group being non-zero
then we have $\Ext^{2(a+b)}(M,Q)$ is non-zero and so $\gfd(Q)=2(a+b)$.

Now
$$\Ext^{4(a+b)+2}\bigl(\Delta_p(\mu),\nabla_p(\mu)\bigr)
\cong \Ext^{4(a+b)}(Q^\circ, Q)$$
as $\gfd(M_\mu)=a < 2(a+b)-1$ for $b \ge 1$.
We let $\lambda=p(a+1,b-1)+(p-r-s-3,r)$.
We also  have
$$
\Ext^{4(a+b)}\bigl(Q^\circ,\nabla_p(\lambda)\bigr)
\cong 
\Ext^{4(a+b)}\bigl(\Delta_p(\lambda),Q\bigr).
$$
Now
$ \Ext^{4(a+b)}(\Delta_p(\lambda),\nabla_p(\lambda)\bigr)\not \cong 0$
by induction.
Hence 
$$
\Ext^{4(a+b)}\bigl(\Delta_p(\lambda),Q\bigr)
\cong
\Ext^{4(a+b)}\bigl(Q^\circ,\nabla_p(\lambda)\bigr)\not\cong0$$
and so $\Ext^{4(a+b)}(Q^\circ, Q)$ is non-zero, using $M=Q^\circ$ in
the sequence above. 
This implies that $\gfd\bigl(\nabla_p(\mu)\bigr)=\gfd(Q)+1= 2(a+b)+1$
and $\Ext^{4(a+b)+2}\bigl(\Delta_p(\mu),\nabla_p(\mu)\bigr)$ is non-zero
as required.


If $b=0$ then 
$$Q\cong \bigl(\nabla(a-1,0)\frob \otimes L(p-r-s-3,r)\bigr)
\oplus \bigl( \nabla(a,0)\frob \otimes L(r,s)\bigr).$$
Hence $\gfd(Q)=2a$ by induction and
$\gfd\bigl(\nabla_p(\mu)\bigr)=2a+1$.
Also
$$\begin{array}{l}\Ext^{4a +2}\bigl(\Delta_p(\mu),\nabla_p(\mu)\bigr)
\cong \Ext^{4a}(Q^\circ, Q)
\\ \hspace{10pt}
\cong \Ext^{4a}\bigl(\Delta(a,0)\frob \otimes L(r,s),\nabla(a,0)\frob
\otimes L(r,s)\bigr) 
\end{array}
$$
which is non-zero by induction.

In all cases for both $\lambda$ and $\mu$ as defined above
we have by  Lemma~\ref{lem:ind2} that 
$\gfd\bigl(L(\lambda)\bigr)=\gfd(\nabla_p(\lambda)\bigr)$
and $\gfd\bigl(L(\mu)\bigr)=\gfd(\nabla_p(\mu)\bigr)$.
This completes the induction. 
\end{proof}


We now consider the case where $\lambda$ lies on a wall but is not a
Steinberg weight.
\begin{lem}\label{lem:mres}
Let $\lambda = p(a,b)+(s,p-1)$ with $(a,b)\in X^+$ and
$0\le s\le p-2$. We define $M_\lambda$ to be the (unique up to
equivalence) non-split extension 
$$0\rarr \nabla(a,b)\frob\otimes L(s,p-1)\rarr M_{\lambda} \rarr
\nabla(a,b-1)\frob\otimes L(r,s) \rarr 0$$
with $r+s=p-2$.
Then $M_\lambda$ has good resolution
$$0\rarr M_\lambda \rarr \nabla(\lambda_0) \rarr \nabla(\lambda_1) 
\rarr \cdots \rarr \nabla(\lambda_a)\rarr 0$$
where 
$ \lambda_{2i}=p(a-2i,b+i) +(s,p-1)$ and 
$\lambda_{2i+1}=p(a-2i-1,b+i+1) +(r,s)$
for $i$ an integer between $0$ and $\lfloor \frac{a}{2}\rfloor$. 

Similarly, let 
$\mu = p(a,b)+(r,s)$ with $(a,b)\in X^+$ and
$r+ s = p-2$. We define $M_\mu$ to be the (unique up to equivalence) 
non-split extension 
$$0\rarr \nabla(a,b)\frob\otimes L(r,s)\rarr M_{\mu} \rarr
\nabla(a,b-1)\frob\otimes L(p-1,r) \rarr 0.$$
Then $M_\mu$ has good resolution
$$0\rarr M_\mu \rarr \nabla(\mu_0) \rarr \nabla(\mu_1) \rarr \cdots \rarr
\nabla(\mu_a)\rarr 0$$
where 
$ \mu_{2i}=p(a-2i,b+i) +(r,s)$ and 
$\mu_{2i+1}=p(a-2i-1,b+i) +(p-1,r)$
for $i$ an integer between $0$ and $\lfloor \frac{a}{2}\rfloor$. 
\end{lem}
\begin{proof}
Similar to that of Lemma~\ref{lem:mzero}
\end{proof}

\begin{cor}\label{cor:gfdmres}
For $M_{\lambda}$ and $M_{\mu}$ as defined above we have
$\gfd(M_{\lambda})=\gfd(M_{\mu})=a$.
Furthermore for $\tau \in X^+$ we have
$$\Ext^a\bigl(\Delta(\tau), M_\lambda \bigr) \cong 
\Hom \bigl(\Delta(\tau), \nabla(\lambda_a)\bigr)\cong \delta_{\tau
\lambda_a}k$$
and
$$\Ext^a\bigl(\Delta(\tau), M_\mu \bigr) \cong 
\Hom \bigl(\Delta(\tau), \nabla(\mu_a)\bigr)\cong \delta_{\tau
\mu_a}k.$$
\end{cor}

The next Proposition shows that $\gfd\bigl(L(\lambda)\bigr)
=g(\lambda)$ for $\lambda$ lying on a wall.

\begin{thm}\label{thm:gfdwall}
Let $\nu$ be a non-Steinberg weight on a wall. (So
$\nu$ is in a primitive block.) Then we have
$$ \gfd\bigl(L(p(a,b)+\nu)\bigr)
=\gfd\bigl(\nabla(a,b)\frob\otimes L(\nu)\bigr)=a+b.
$$
Furthermore
$$\Ext^{2(a+b)}\bigl(\Delta(a,b)\frob\otimes L(\nu),
\nabla(a,b)\frob\otimes L(\nu)\bigr) \not \cong 0.$$
\end{thm}
\begin{proof} We proceed by induction on $a+b$.
For $a+b=0$ we have $L(\nu)=\nabla(\nu)=\Delta(\nu)$ and 
$\gfd\bigl(L(\nu)\bigr)=\gfd\bigl(\nabla(\nu)\bigr)=0.$ 

For $a+b=1$, consider $a=1$, $b=0$ (the other case being similar).

Case (i): $\nu=(r,s)$ with $r+s=p-2$. Here we have
$L(p+r,s)\cong \nabla(1,0)\frob
\otimes L(r,s)$. We have a non-split short exact sequence~\cite[lemma
3.2.4 (iv)]{yehia}
$$0\rarr L(p+r,s)\rarr \nabla(p+r,s)\rarr \nabla(s,p-1)\rarr 0,$$
so $\gfd\bigl(L(p+r,s)\bigr)=1$ and
$\Ext^2\bigl(L(p+r,s),L(p+r,s)\bigr)\not\cong 0 $ by
Lemma~\ref{lem:globtwo}.
 
Case (ii): $\nu=(p-1,r)$. This case follows exactly as in Case (i).

Case (iii): $\nu=(s,p-1)$. Here 
$L(p+s,p-1)\cong \nabla(1,0)\frob\otimes L(p-1,r)$. We have a non-split short exact
sequence~\cite[lemma 3.2.4 (ii) and (vii)]{yehia}
$$0\rarr L(p+s,p-1)\rarr \nabla(p+s,p-1) \rarr \nabla(r,p+s)\rarr 0$$
so $\gfd\bigl(L(p+s,p-1)\bigr)= 1$ and  
$\Ext^2\bigl(L(p+s,p-1),L(p+s,p-1)\bigr)\not\cong 0$ by
Lemma~\ref{lem:globtwo}.

Now suppose that $a+b\ge 2$.
We let $\lambda= p(a,b)+\nu$. We will only consider the case with
$\nu=(r,s)$ or $(s,p-1)$. The other case with $\nu=(p-1,r)$ is exactly
dual to the case of $\nu=(s,p-1)$.
We define $\mu$ by
$$\mu = \left\{\begin{array}{ll}
                p(a,b-1) +(p-1,r) &\mbox{if $\nu=(r,s)$}\\
                p(a,b-1) +(r,s)   &\mbox{if $\nu=(s,p-1)$}
                \end{array}
        \right.
$$
where $r+s$ equals $p-2$.
We have a short exact sequence
\begin{equation}\label{Mses3}
0\rarr \nabla_p(\lambda) \rarr M_\lambda \rarr \nabla_p(\mu)\rarr 0
\end{equation}
with $M_\lambda$ as defined in Lemma~\ref{lem:mres}.

Case (i): $b=0$. Here $\nabla_p(\mu)=0$ so $\nabla_p(\lambda)=M_\lambda$.
Corollary~\ref{cor:gfdmres} gives 
$\gfd(\nabla_p(\lambda))=a$ and 
$$\Ext^{2a}\bigl(\Delta_p(\lambda),
\nabla_p(\lambda)\bigr) \cong \Hom\bigl(\Delta(\lambda_a),
\nabla(\lambda_a)\bigr)\cong k$$
where $\lambda_a$ is defined as in Lemma~\ref{lem:mres} and using the
good resolution for $M_\lambda$ and its $^\circ$-dual.

Case (ii): $b=1$.
We know by Case (i) that 
$\gfd\bigl(\nabla_p(\mu)\bigr)=a$ 
and Corollary~\ref{cor:gfdmres} give $\gfd(M_\lambda)=a$. 
Hence
Lemma~\ref{lem:gandi} applied to sequence~(\ref{Mses3}) gives 
$\gfd\bigl(\nabla_p(\lambda)\bigr) \le a+1$.
Using Case (i) and Corollary~\ref{cor:gfdmres}
we have the commutative diagram
$$\xymatrix@R=10pt@C=10pt{
{\Ext^{a}\bigl(\Delta(\tau),M_{\lambda}\bigr)}\ar@{->}[r] \ar@{=}[d] 
&{\Ext^{a}\bigl(\Delta(\tau),\nabla_p(\mu)\bigr)}
\ar@{->}[r] \ar@{=}[d] 
&{\Ext^{a+1}\bigl(\Delta(\tau),\nabla_p(\lambda) \bigr)}\ar@{->}[r] \ar@{=}[d] 
& 0 \ar@{=}[d] \\
{\delta_{\tau\lambda_{a}}k}\ar@{->}[r]  
& {\delta_{\tau\mu_{a-1}}k}\ar@{->}[r]  
&{\Ext^{a+1}\bigl(\Delta(\tau),\nabla_p(\lambda) \bigr)}\ar@{->}[r]  
& 0 
}$$
But $\lambda_{a}\ne \mu_a$ and hence we have
$\gfd\bigl(\nabla_p(\lambda)\bigr)= a+1$.

We now wish to show that
$\Ext^{2a+2}\bigl(\Delta_p(\lambda),\nabla_p(\lambda)\bigr)$ is
non-zero. We have just shown that
$\wfd\bigl(\Delta_p(\lambda)\bigr)=\gfd\bigl(\nabla_p(\lambda)\bigr)=a+1$.
Also by Case (i) we have
$\wfd\bigl(\Delta_p(\mu)\bigr)=\gfd\bigl(\nabla_p(\mu)\bigr)= a$.
Hence, using Lemma~\ref{lem:gandi}, the long exact sequence from 
sequence~(\ref{Mses3}) gives us
\begin{equation}\label{longses}
\Ext^{2a+1}\bigl(\Delta_p(\lambda),M_\lambda\bigr) 
\rarr \Ext^{2a+1}\bigl(\Delta_p(\lambda), \nabla_p(\mu)\bigr)
\rarr \Ext^{2a+2}\bigl(\Delta_p(\lambda), \nabla_p(\lambda)\bigr)
\rarr 0
\end{equation}
Now $\Ext^{2a+1}\bigl(\Delta_p(\lambda), M_{\lambda}\bigr)
\cong \Ext^{a+2}\bigl(\Delta_p(\lambda),M_{\lambda_{a-1}}\bigr)$
using the good resolution for $M_\lambda$ from\break Lemma~\ref{lem:mres}.
We have the following short exact sequence for $M_{\lambda_{a-1}}$
$$0\rarr M_{\lambda_{a-1}} \rarr \nabla(\lambda_{a-1}) \rarr
\nabla(\lambda_{a}) \rarr 0$$
so the long exact sequence gives us
$$\Ext^{a+1}\bigl(\Delta_p(\lambda), \nabla(\lambda_a)\bigr) \rarr
\Ext^{a+2}\bigl(\Delta_p(\lambda), M_{\lambda_{a-1}}\bigr) \rarr 0$$
as $\wfd(\Delta_p(\lambda))=a+1$.
But $\Ext^{a+1}\bigl(\Delta_p(\lambda), \nabla(\lambda_a)\bigr)\cong
\delta_{\lambda_a \mu_a}k=0$. Hence $\Ext^{a+2}\bigl(\Delta_p(\lambda),
M_{\lambda_{a-1}}\bigr) \cong 0$. Sequence~(\ref{longses}) then gives
us 
$$\begin{array}{rl}
\Ext^{2a+2}\bigl(\Delta_p(\lambda), \nabla_p(\lambda)\bigr) 
&\cong
\Ext^{2a+1}\bigl(\Delta_p(\lambda), \nabla_p(\mu)\bigr)
\\ &\cong 
\Ext^{2a+1}\bigl(\Delta_p(\mu), \nabla_p(\lambda)\bigr) .
\end{array}$$


We now need to show that
$$
\Ext^{2a+1}\bigl(\Delta_p(\mu), \nabla_p(\lambda)\bigr) 
\cong
\Ext^{2a}\bigl(\Delta_p(\mu), \nabla_p(\mu)\bigr)
$$
where the last $\Ext$ group is non-zero by induction.
But we may repeat the argument above in one less degree using $\Delta_p(\mu)$
in place of $\Delta_p(\lambda)$ to get the required isomorphism.

Case (iii): $b\ge 2$.
Here $\gfd\bigl(\nabla_p(\mu)\bigr)=a+b-1> \gfd(M)=a$ so
Lemma~\ref{lem:gandi}
gives us $\gfd\bigl(\nabla_p(\lambda)\bigr)=a+b$,
and 
$$\begin{array}{rl}
\Ext^{2(a+b)}\bigl(\Delta_p(\lambda),\nabla_p(\lambda)\bigr)
&\cong \Ext^{2(a+b-1)}\bigl(\Delta_p(\mu),\nabla_p(\lambda)\bigr)
\\ &\cong \Ext^{2(a+b)-2}\bigl(\Delta_p(\mu), \nabla_p(\mu)\bigr)
\not \cong 0
\end{array}$$
by induction.

Lemma~\ref{lem:ind2} then completes the proof.
\end{proof}

We now consider the case where $\lambda$ is not primitive.
\begin{cor}\label{cor:gfdst}
Suppose $\lambda$ is a dominant weight and $\lambda=p^d \lambda_1+
(p^d-1,p^d-1)$ for some $d \in \N$ and $\lambda_1\in X^+$. 
Then $\gfd\bigl(L(\lambda)\bigr)= \gfd\bigl(L(\lambda_1)\bigr)$ and 
$$\Ext^{2 \gfd(L(\lambda))}\bigl(L(\lambda), L(\lambda)\bigr) \not \cong 0.$$
\end{cor}
\begin{proof}
Suppose $\lambda=(a,b)$, and $\lambda_1=(a_1,b_1)$.
We have by~\cite[section 4, theorem]{donk4} that 
$\B(a+b,b,0)$ is Morita equivalent to $\B(a_1+b_1,b_1,0)$ in $S(2,r_1)$
with $r_1= a_1+2b_1$.  But $\lambda_1$ is
primitive and the result follows by Theorems~\ref{thm:gfdalc} and
~\ref{thm:gfdwall}.
\end{proof}


\begin{cor}\label{cor:s3gfd}
Given $(a_1,a_2,a_3)\in \partn(3,r)$,
 we let $(a_1-a_2,a_2-a_3)
= p^{d}\lambda+  (p^d-1,p^d-1)$
with $\lambda\in X^+$ and $d \in \N$.
We also let $L(a_1,a_2,a_3)$ 
be the irreducible
module of highest weight $(a_1,a_2,a_3)$ for $S(3,r)$.
Then 
$$\gfd\bigl(L(a_1,a_2,a_3)\bigr)= g(\lambda).$$
Moreover
$$\Ext^{2 g(\lambda)}_{S(3,r)}\bigl(L(a_1,a_2,a_3),L(a_1,a_2,a_3)\bigr)
\not\cong 0.$$
\end{cor}
\begin{proof}
We let $g=g(\lambda)$.
Now Theorems~\ref{thm:gfdalc} and~\ref{thm:gfdwall} and
Corollary~\ref{cor:gfdst}
 give us
$$\begin{array}{l}
\Ext_{S(3,r)}^i\bigl(\Delta(b_1,b_2,b_3), L(a_1,a_2,a_3)\bigr)
\\
\cong\Ext_{G}^i\bigl(\Delta(b_1-b_2,b_2-b_3), L(a_1-a_2,
a_2-a_3)\bigr)
\cong 0
\end{array}
$$
if $i> g$ and so we have $\gfd\bigl(L(a_1,a_2,a_3)\bigr) \le g$.
We also have 
$$
\begin{array}{l}
\Ext^{2g}_{S(3,r)}\bigl(L(a_1,a_2,a_3),L(a_1,a_2,a_3)\bigr)
\\\cong\Ext^{2g}_G\bigl(L(a_1-a_2,a_2-a_3),L(a_1-a_2,a_2-a_3)\bigr)
\not\cong 0.
\end{array}
$$
Now Lemma~\ref{lem:wandg} gives 
$\wfd(L(a_1,a_2,a_3))+\gfd(L(a_1,a_2,a_3))\ge 2 g$.
But $\wfd(L(a_1,a_2,a_3))= \gfd(L(a_1,a_2,a_3))$, and 
so we have $\gfd(L(a_1,a_2,a_3))= g$, as required.
\end{proof}

\begin{thm}\label{thm:s3glob}
The global dimension of $S(3,r)$ is twice its good filtration
dimension and  is given as follows
$$\begin{array}{llrcl}
\mbox{for}& p=2 & \glob(S(3,r)) &= & 2\lfloor\frac{r}{2}\rfloor \\
 & & & & \\
\mbox{for}& p=3 & \glob(S(3,r)) &= & \left\{ \begin{array}{ll}
                                4\!\left(\frac{r}{3}\right) &\mbox{if $r
                                        \equiv 0\pmod 3$} \\
                                2\lfloor\frac{r}{3}\rfloor &\mbox{if $r
                                        \not\equiv 0\pmod 3$}
                                \end{array}\right. \\
 & & & & \\
\mbox{for}& p\ge5 & \glob(S(3,r)) &= &4\lfloor\frac{r}{p}\rfloor.
\end{array}$$
\end{thm}
\begin{proof}
Using Corollary~\ref{cor:s3gfd},
the same argument as in the second paragraph of the proof
of Theorem~\ref{thm:s2glob} gives us $\glob(S(3,r)) =2 \gfd(S(3,r))$.

For $\SL_3$, $(r,0)$ is always primitive.
Lemmas~\ref{lem:gdecr} and~\ref{lem:wandg} mean we need only 
consider weights in $S(3,r)$
which are maximal in their block in determining $\gfd(S(3,r))$.

Case (i): $p=2$. Here all the corresponding $\SL_3$ weights lie on
walls. Also all the irreducible modules corresponding to weights in
$\partn(3,r)$ have smaller or equal good filtration dimension to
$(r,0)$, so $\glob(S(3,r))= 2\gfd(L(r,0,0))= 2 \lfloor \frac{r}{2} \rfloor$.

Case (ii): $p=3$. If $r \equiv 0\pmod 3$ then $(r,0)$ lies inside a
lower alcove and $\glob(S(3,r))= 4\left(\frac{r}{3}\right)$.
If $r \equiv 1$ or $2\pmod3$ then $(r,0)$ lies on a wall.
We claim all the $\SL_3$ weights corresponding to partitions 
in $\partn(3,r)$ lie on a wall and then all irreducible modules
corresponding to weights in $\partn(3,r)$ will have smaller
or equal good filtration dimension.
If $(a,b)$ lies inside an alcove then $(a,b) \equiv (0,0)$ or
$(2,2)\pmod 3$. So in both cases $a+2b\equiv 0\pmod p$. 
For $(a,b)$ to come from a partition in $\partn(3,r)$
we need  $(a+b+c,b+c,c) \in \partn(3,r)$ for some $c
\in \N$ with $a+2b+3c=r$. 
But we have 
$a+2b\equiv 0\pmod 3$ so $r \equiv 0 \pmod 3$,
and this is a contradiction to our assumption on $r$.
Thus all the $\SL_3$ weights corresponding to partitions in 
$\partn(3,r)$ lie on a wall and hence
$\glob(S(3,r))= 2\lfloor\frac{r}{3}\rfloor$.

Case (iii): $p \ge 5$.
If $(r,0)$ is not on a wall then $\glob(S(3,r))=
4\lfloor\frac{r}{p}\rfloor$.
If $(r,0)$ is on a wall then $r \equiv -1$ or $-2\pmod p$.
But then $\partn(3,r)$ contains $(r-2,1,1)$ and this corresponds to
the $\SL_3$ weight $(r-3,0)$ which lies inside an
alcove.
Hence $\glob(S(3,r))= 4\lfloor\frac{r-3}{p}\rfloor= 4\lfloor\frac{r}{p}\rfloor$.
\end{proof}

We now consider the quantum case with $n=3$ and $l\ge 2$. 
The cohomological theory of quantum groups and their $q$-Schur
algebras can be found in~\cite{donkquant}.
We need the appropriate generalisation of the
$p$-filtration of $\nabla(\lambda)$ in Theorem~\ref{thm:pfiltra},
for which we need the generalisation of Proposition~\ref{propn:yehia}. 
We expect that this would replace $p$ by $l$ in all cases.
Once we have a $G_1B$ composition series of $\Zone(\lambda)$ 
then~\cite[proposition 5.2]{cox3} would give us the required $l$-filtration of
$\nabla(\lambda)$. This would then give a quantum version of the
$M_\lambda$'s used extensively in Section 5 and then all the Theorems
in this section would generalise. 
We expect that the $G_1B$ composition series of $\Zone(\lambda)$ will
have the same weights in terms of relative alcoves as in the
classical case.
This would then give us the quantum version of Theorem~\ref{thm:s3glob}
with $p$ replaced by $l$ and $5$ replaced by $4$.


\section*{Acknowledgments}
I would like to thank my PhD supervisor, Stephen Donkin for his great help and
encouragement. I would also like to thank Anton Cox for many helpful
suggestions and comments.
This research was supported by the Association of Commonwealth
Universities and the British Council.

\providecommand{\bysame}{\leavevmode\hbox to3em{\hrulefill}\thinspace}
\providecommand{\MR}{\relax\ifhmode\unskip\space\fi MR }
\providecommand{\MRhref}[2]{%
  \href{http://www.ams.org/mathscinet-getitem?mr=#1}{#2}
}
\providecommand{\href}[2]{#2}

%

\end{document}